\documentclass[12pt]{article}
\usepackage{amssymb, amsmath, latexsym, a4}
\usepackage[latin1]{inputenc}
\usepackage[all]{xy}
 \usepackage[T1]{fontenc}

\begin{document}

\newcommand{\rdg}{\hfill $\Box $}

\newtheorem{De}{Definition}[section]
\newtheorem{Th}[De]{Theorem}
\newtheorem{Pro}[De]{Proposition}
\newtheorem{Le}[De]{Lemma}
\newtheorem{Co}[De]{Corollary}
\newtheorem{Rem}[De]{Remark}
\newtheorem{Ex}[De]{Example}
\newtheorem{Exo}[De]{Exercises}
\newcommand{\tp}{\otimes}
\newcommand{\N}{\mathbb{N}}
\newcommand{\Z}{\mathbb{Z}}
\newcommand{\K}{\mathbb{K}}
\newcommand{\op}{\oplus}
\newcommand{\n}{\underline n}
\newcommand{\es}{{\frak S}}
\newcommand{\ef}{\frak F}
\newcommand{\qu}{\frak Q} \newcommand{\ga}{\frak g}
\newcommand{\la}{\lambda}
\newcommand{\ig}{\frak Y}
\newcommand{\te}{\frak T}
\newcommand{\cok}{{\sf Coker}}
\newcommand{\Hom}{{\sf Hom}}
\newcommand{\im}{{\sf Im}}
\newcommand{\ext}{{\sf Ext}}
\newcommand{\ho}{{\sf H_{AWB}}}
\newcommand{\HH}{{\sf Hoch}}
\newcommand{\adu}{{\rm AWB}^!}

\newcommand{\ele}{\cal L} \newcommand{\as}{\cal A} \newcommand{\ka}{\cal
K}\newcommand{\eme}{\cal M} \newcommand{\pe}{\cal P}

\newcommand{\pn}{\par \noindent}
\newcommand{\pbn}{\par \bigskip \noindent}
\bigskip\bigskip

\centerline {\Large {\bf  On the universal $\alpha$-central extension of the}}

\centerline {\Large {\bf  semi-direct product of Hom-Leibniz algebras}}

\

\centerline {\bf J. M. Casas$^{(1)}$ and N. Pacheco Rego$^{(2)}$}

\bigskip \bigskip
\centerline{$^{(1)}$Dpto.  Matem\'atica Aplicada I, Univ. de Vigo, 36005 Pontevedra, Spain}
\centerline{e-mail address: jmcasas@uvigo.es}

\centerline{$^{(2)}$IPCA, Dpto. de Ciências, Campus do IPCA,
 Lugar do Aldão}
\centerline{4750-810 Vila Frescainha, S. Martinho, Barcelos,
 Portugal}
\centerline{e-mail address: nrego@ipca.pt}

\bigskip \bigskip \bigskip \bigskip

\par
{\bf Abstract} We introduce Hom-actions, semidirect product and establish the equivalence between split extensions and the semi-direct product extension of Hom-Leibniz algebras. We analyze the functorial properties of the universal ($\alpha$)-central extensions of ($\alpha$)-perfect Hom-Leibniz algebras. We establish under what conditions an automorphism or a derivation can be lifted in an $\alpha$-cover and we analyze  the universal $\alpha$-central extension of the semi-direct product of two $\alpha$-perfect Hom-Leibniz algebras.

\bigskip \bigskip


\bigskip \bigskip
{\it A. M. S. Subject Class. (2010):} 17A32, 16E40, 17A30

\section{Introduction}

Hom-Lie algebras were introduced in \cite{HLS} as Lie algebras whose Jacobi identity is twisted by means of a map. This fact occurs in different applications in models of quantum phenomena or in analysis of complex systems and processes exhibiting complete or partial scaling invariance.

From the introducing paper, the investigation of several kinds of Hom-structures is in progress (for instance, see \cite{AM, AMM, AMS,AMS1, Ma, MS, Yau4, Yu} and references given therein). Naturally, the non-skew-symmetric version of Hom-Lie algebras, the so called Hom-Leibniz algebras, was considered as well (see \cite{AMM, CIP, Cheng, Is, MS2, MS, Yau1}). A Hom-Leibniz algebra is un triple $(L,[-,-],\alpha_L)$ consisting of a $\mathbb{K}$-vector space $L$, a bilinear map $[-,-] : L \times L \to L$ and a  homomorphism of
$\mathbb{K}$-vector spaces $\alpha_L : L \to L$ satisfying the Hom-Leibniz identity:
$$[\alpha_L(x),[y,z]]=[[x,y], \alpha_L(z)]-[[x,z], \alpha_L(y)]$$
for all $x, y, z \in L$. When $\alpha_L=Id$, the definition of Leibniz algebra \cite{Lo} is recovered. If the bracket is skew-symmetric, then we recover the definition of Hom-Lie algebra \cite{HLS}.

The goal of the present paper is to continue with the investigations on universal ($\alpha$)-central extensions of ($\alpha$)-perfect Hom-Leibniz algebras initiated in \cite{CIP}.
In concrete, we consider the extension of  results about the universal central extension of the semi-direct product of Leibniz algebra in \cite{CC} to the framework of Hom-Leibniz algebras.

To do so, we organize the paper as follows: an initial section recalling the background material on Hom-Leibniz algebras. We introduce the concepts of Hom-action and semi-direct product and we prove a new result (Lemma \ref{ext rota}) that establishes the equivalence between split extensions and the semi-direct product extension. Section 3 is devoted to analyze the functorial properties of the universal ($\alpha$)-central extensions of ($\alpha$)-perfect Hom-Leibniz algebras. In section 4 we establish under what conditions an automorphism or a derivation can be lifted in an $\alpha$-cover (a central extension $f:(L', \alpha_{L'}) \to (L, \alpha_L)$ where $(L', \alpha_{L'})$ is $\alpha$-perfect ($L' = [\alpha_{L'}(L'), \alpha_{L'}(L')]$)). Final section is devoted to analyze the relationships between the universal $\alpha$-central extension of the semi-direct product of two $\alpha$-perfect Hom-Leibniz algebras, such that one of them Hom-act over the other one, and the semi-direct product of the universal $\alpha$-central extensions of both of them.

\section{Preliminaries on Hom-Leibniz algebras}

In this section we introduce necessary material on  Hom-Leibniz algebras which will be used in subsequent sections.

\begin{De}\label{HomLeib} \cite{MS}
A Hom-Leibniz algebra is a triple $(L,[-,-],\alpha_L)$ consisting of a $\mathbb{K}$-vector space $L$, a bilinear map $[-,-] : L \times L \to L$ and a $\mathbb{K}$-linear map $\alpha_L : L \to L$ satisfying:
\begin{equation} \label{def}
 [\alpha_L(x),[y,z]]=[[x,y],\alpha_L(z)]-[[x,z],\alpha_L(y)] \ \ ({\rm Hom-Leibniz\ identity})
\end{equation}
for all $x, y, z \in L$.

A Hom-Leibniz algebra $(L,[-,-],\alpha_L)$  is said to be multiplicative \cite{Yau1} if the $\mathbb{K}$-linear map $\alpha_L$ preserves the bracket, that is, if $\alpha_L [x,y] = [\alpha_L(x),\alpha_L(y)]$, for all $x, y \in L$.
\end{De}

\begin{Ex}\label{ejemplo 1} \
\begin{enumerate}
\item[a)] Taking $\alpha = Id$ in Definition \ref{HomLeib} we obtain the definition of  Leibniz algebra \cite{Lo}. Hence Hom-Leibniz algebras include Leibniz algebras as a full subcategory, thereby motivating the name "Hom-Leibniz algebras" as a deformation of Leibniz algebras twisted by a homomorphism. Moreover it  is a multiplicative Hom-Leibniz algebra.

    \item[b)] Hom-Lie algebras \cite{HLS} are Hom-Leibniz algebras whose bracket satisfies the condition $[x,x]=0$, for all $x$. So Hom-Lie algebras can be considered as a full subcategory of Hom-Leibniz algebras category. For any multiplicative Hom-Leibniz algebra $(L,[-,-],\alpha_L)$ there is associated the Hom-Lie algebra $(L_{\rm Lie},[-,-],\widetilde{\alpha})$, where $L_{\rm Lie} = L/L^{\rm ann}$, the bracket is the canonical bracket induced on the quotient and  $\widetilde{\alpha}$ is the homomorphism naturally induced by $\alpha$. Here $L^{\rm ann} = \langle \{[x,x] : x \in L \} \rangle$.

    \item[b)] Let $(D,\dashv,\vdash, \alpha_D)$ be a Hom-dialgebra. Then $(D,\dashv,\vdash, \alpha_D)$ is a  Hom-Leibniz algebra with respect to the bracket  $[x,y]=x \dashv y - y \vdash x$, for all $x,y \in A$ \cite{Yau}.

  \item[c)] Let $(L,[-,-])$ be a Leibniz algebra and $\alpha_L:L \to L$  a Leibniz algebra endomorphism. Define $[-,-]_{\alpha} : L \otimes L \to L$ by $[x,y]_{\alpha} = [\alpha(x),\alpha(y)]$, for all $x, y \in L$. Then $(L,[-,-]_{\alpha}, \alpha_L)$ is a multiplicative Hom-Leibniz algebra.

\item[d)] Abelian or commutative Hom-Leibniz algebras are $\mathbb{K}$-vector spaces $L$ with trivial bracket and any linear map $\alpha_L :L \to L$.

\end{enumerate}
\end{Ex}

\begin{De}\label{homo}
A homomorphism of  Hom-Leibniz algebras $f:(L,[-,-],\alpha_L) \to (L',[-,-]',\alpha_{L'})$ is a $\mathbb{K}$-linear map $f : L \to L'$ such that
\begin{enumerate}
\item[a)] $f([x,y]) =[f(x),f(y)]'$
\item [b)] $f \cdot \alpha_L(x) = \alpha_{L'} \cdot  f(x)$
\end{enumerate}
for all $x, y \in L$.

A homomorphism of multiplicative Hom-Leibniz algebras is a homomorphism of the underlying Hom-Leibniz algebras.
\end{De}
\medskip

In the sequel we refer Hom-Leibniz algebra to a multiplicative Hom-Leibniz algebra and we shall use the shortened notation $(L,\alpha_L)$ when there is not confusion with the bracket operation.

\begin{De}
Let $(L,[-,-],\alpha_L)$ be a Hom-Leibniz algebra. A  Hom-Leibniz subalgebra $(H, \alpha_H)$ is a linear subspace $H$ of $L$, which is closed for the bracket and invariant by $\alpha_L$, that is,
\begin{enumerate}
\item [a)] $[x,y] \in H,$ for all $x, y \in H$,
\item [b)] $\alpha_L(x) \in H$, for all $x \in H$ ($\alpha_H = \alpha_{L \mid}$).
\end{enumerate}

A  Hom-Leibniz subalgebra $(H, \alpha_H)$ of $(L, \alpha_L)$ is said to be a  two-sided Hom-ideal if $[x,y], [y,x] \in H$, for all $x \in H, y \in L$.

If $(H,\alpha_H)$ is a two-sided  Hom-ideal of $(L,\alpha_L)$, then the quotient $L/H$ naturally inherits a structure of Hom-Leibniz algebra with respect to the endomorphism $\widetilde{\alpha} : L/H  \to L/H, \widetilde{\alpha}(\overline{l}) =\overline{\alpha_L(l)}$, which is said to be the quotient Hom-Leibniz algebra.
\end{De}

So we have defined the category ${\sf Hom-Leib}$ (respectively, ${\sf Hom-Leib_{\rm mult}})$ whose objects are Hom-Leibniz (respectively, multiplicative Hom-Leibniz) algebras and whose morphisms are the homomorphisms of Hom-Leibniz (respectively, multiplicative Hom-Leibniz) algebras.
There is an obvious inclusion functor $inc : {\sf Hom-Leib_{\rm mult}} \to {\sf Hom-Leib}$. This functor has as left adjoint the multiplicative functor $(-)_{\rm mult} : {\sf Hom-Leib} \to {\sf Hom-Leib_{\rm mult}}$ which assigns to a Hom-Leibniz algebra $(L,[-,-],\alpha_L)$ the multiplicative Hom-Leibniz algebra $(L/I,[-,-],\tilde{\alpha})$, where $I$ is the two-sided ideal of $L$ spanned by the elements $\alpha_L[x,y]-[\alpha_L(x),\alpha_L(y)]$, for all $x, y \in L$.

\begin{De}
Let $(H, \alpha_H)$ and $(K, \alpha_K)$ be two-sided Hom-ideals of a Hom-Leibniz algebra $(L,[-,-],\alpha_L)$. The commutator  of $(H, \alpha_H)$ and $(K, \alpha_K)$, denoted by $([H,K],\alpha_{[H,K]})$, is the Hom-Leibniz subalgebra of $(L,\alpha_L)$ spanned by the brackets $[h,k], h \in H, k \in K$.
\end{De}

Obviously, $[H,K] \subseteq H \cap K$ and $[K,H] \subseteq H \cap K$. When $H = K =L$, we obtain the definition of derived Hom-Leibniz subalgebra. Let us observe that, in general, $([H,K],\alpha_{[H,K]})$ is not   a Hom-ideal, but if $H, K \subseteq \alpha_L(L)$, then $([H,K],\alpha_{[H,K]})$ is a two-sided ideal of $(\alpha_L(L), \alpha_{L \mid})$. When $\alpha = Id$, the classical notions are recovered.

\begin{De}
Let $(L,[-,-],\alpha_L)$ be a Hom-Leibnz algebra. The subspace $Z(L) = \{ x \in L \mid [x, y] =0 = [y,x], \text{for\ all}\ y \in L \}$ is said to be the center of $(L,[-,-],\alpha_L)$.

When $\alpha_L : L \to L$ is a surjective homomorphism, then  $Z(L)$ is a Hom-ideal of $L$.
\end{De}


\subsection{Hom-Leibniz actions}
\begin{De} \label{Hom accion} Let $\left(  L ,\alpha_{L}\right)$ and $\left(  M,\alpha_{M}\right)$ be Hom-Leibniz algebras. A (right) Hom-action of $\left(  L, \alpha_{L}\right)$ over $\left(  M, \alpha_{M}\right)$ consists of two bilinear maps,
$\lambda:L\otimes M\to M$, $\lambda\left(  l\otimes m\right)  =l \centerdot  m$, and $\rho:M\otimes L\to M$, $\rho \left(  m\otimes
l\right)  =m \centerdot  l$, satisfying the following identities:
\begin{enumerate}
\item[a)] $\alpha_{M}\left(  m\right)  \centerdot  \left[  x,y\right]  =\left(
m \centerdot  x\right)  \centerdot \alpha_{L}\left(  y\right)  -\left(  m \centerdot  y\right)  \centerdot  \alpha_{L}\left(  x\right);$

\item[b)] $\alpha_{L}\left(  x\right)  \centerdot  \left(  m \centerdot y\right)  =\left(
x \centerdot m\right)  \centerdot  \alpha_{L}\left(  y\right)  -\left[  x,y\right]  \centerdot  \alpha
_{M}\left(  m\right);$

\item[c)] $\alpha_{L}\left(  x\right) \centerdot  \left(  y \centerdot  m\right)  =\left[
x,y\right]  \centerdot  \alpha_{M}\left(  m\right)  -\left(  x \centerdot  m\right)  \centerdot  \alpha
_{L}\left(  y\right);$

\item[d)] $\alpha_{L}\left(  x\right) \centerdot  \left[  m,m^{\prime}\right]  =\left[
 x \centerdot  m,\alpha_{M}\left(  m^{\prime}\right)  \right]  -\left[
 x \centerdot  m^{\prime},\alpha_{M}\left(  m\right)  \right];$

\item[e)] $\left[  \alpha_{M}\left(  m\right)  ,  m^{\prime} \centerdot  x
\right]  =\left[  m,m^{\prime}\right]  \centerdot  \alpha_{L}\left(  x\right)  -\left[
 m \centerdot  x,\alpha_{M}\left(  m^{\prime}\right)  \right];$

\item[f)] $\left[  \alpha_{M}\left(  m\right)  ,  x \centerdot  m^{\prime}
\right]  =\left[    m \centerdot  x,\alpha_{M}\left(  m^{\prime}\right)
\right]  -\left[  m,m^{\prime}\right] \centerdot  \alpha_{L}\left(  x\right);$

\item[g)] $\alpha_{M}\left(  x \centerdot  m\right)  =\alpha_{L}\left(  x\right) \centerdot \alpha_{M}\left(  m\right);$

\item[h)] $\alpha_{M}\left(  m \centerdot  x\right)  =\alpha_{M}\left(  m\right) \centerdot  \alpha_{L}\left(  x\right);$
\end{enumerate}
for all $x,y\in L$ and $m,m^{\prime}\in M$.

When $(M, \alpha_M)$ is an abelian Hom-Leibniz algebra, that is the bracket on $M$ is trivial, then  the Hom-action is called Hom-representation.
\end{De}

\begin{Ex}\label{Hom accion Leib} \
\begin{enumerate}
\item[a)] Let $M$ be a  representation of a  Leibniz algebra L \cite{LP}. Then $(M,Id_M)$ is a  Hom-representation of the  Hom-Leibniz algebra $(L,Id_L)$.

\item[b)]  Let $\left(  K, \alpha_{K}\right)$ be a Hom-Leibniz subalgebra of a  Hom-Leibniz algebra $\left(  L, \alpha_{L}\right)$ (even $\left(  K, \alpha_{K}\right)  =\left(  L, \alpha_{L}\right)  $) and $\left(  H, \alpha_{H}\right)$ a two-sided Hom-ideal of $\left(  L, \alpha_{L}\right)$. There exists a Hom-action of $\left(  K,\alpha_{K}\right)$ over $\left(  H, \alpha_{H}\right)$ given by the bracket in $\left(  L,\alpha_{L}\right).$

\item[c)] An abelian sequence of Hom-Leibniz algebras is an exact sequence of Hom-Leibniz algebras  $0 \to (M,\alpha_M) \stackrel{i}\to (K,\alpha_K) \stackrel{\pi}\to (L,\alpha_L) \to 0$, where $(M,\alpha_M)$ is an abelian Hom-Leibniz algebra, $\alpha_K \cdot i = i \cdot \alpha_M$ and $\pi \cdot \alpha_K = \alpha_L \cdot \pi$.

      An abelian sequence induces a Hom-representation structure of $(L,\alpha_L)$ over $(M,\alpha_M)$ by means of the actions given by $\lambda : L \otimes M \to M, \lambda(l,m)=[k,m], \pi(k)=l$ and $\rho : M \otimes L \to M, \rho(m,l)=[m,k], \pi(k)=l$.
\end{enumerate}
\end{Ex}

\begin{De} \label{producto semidirecto}
Let $\left(  M,\alpha_{M}\right)$ and  $\left(  L,\alpha_{L}\right)$ be Hom-Leibniz algebras together with a Hom-action of $\left(  L,\alpha_{L}\right)$ over $\left(  M,\alpha_{M}\right)$.  Its semi-direct product $\left(  M\rtimes L,\widetilde{\alpha}\right)$ is the Hom-Leibniz algebra with  underlying $\mathbb{K}$-vector space $M\oplus L$,  endomorphism $\widetilde{\alpha}:M\rtimes L\to M\rtimes L$ given by $\widetilde{\alpha} \left(  m,l\right) = \left(  \alpha_{M}\left(m\right)  ,\alpha_{L}\left(  l\right)  \right)$ and bracket
$$\left[  \left(  m_{1},l_{1}\right)  ,\left(  m_{2},l_{2}\right)
\right]  =\left( [m_1,m_2] +  \alpha_{L}\left(  l_{1}\right)  \centerdot  m_{2}+m_{1} \centerdot
\alpha_{L}\left(  l_{2}\right)  ,\left[  l_{1},l_{2}\right]  \right).$$
\end{De}

Let $\left(  M,\alpha_{M}\right)$ and  $\left(  L,\alpha_{L}\right)$ be Hom-Leibniz algebras with a Hom-action of $\left(  L,\alpha_{L}\right)$ over $\left(  M,\alpha_{M}\right)$, then we can construct the  sequence
\begin{equation} \label{extension semidirectoLeib}
0\to \left(  M,\alpha_{M}\right)
\overset{i}{\to}\left(  M\rtimes L,\widetilde{\alpha}\right)
\overset{\pi}{\to}\left(  L,\alpha_{L}\right)  \to 0
\end{equation}
where $i: M \to  M\rtimes L, i(m) =  \left(  m,0\right),$ and $\pi:  M\rtimes L  \to  L,
\pi \left(  m,l\right)   = l$. Moreover, this sequence splits by $\sigma:  L  \to M\rtimes L,
 \sigma(l)  =  \left(  0,l\right),$ that is,  $\sigma$ satisfies $\pi \cdot \sigma=Id_{L}$ and
 $\widetilde{\alpha} \cdot \sigma =\sigma \cdot \alpha_{L}$.

 \begin{De}
Let $(M, \alpha_M)$ and $(L, \alpha_L)$ be Hom-Leibniz algebras such that there is a Hom-action of $(L, \alpha_L)$ over $(M, \alpha_M)$. Two  extensions of $(L, \alpha_L)$ by  $(M, \alpha_M)$,  $0 \to (M, \alpha_M) \stackrel{i} \to (K,\alpha_K) \stackrel{\pi} \to (L, \alpha_L) \to 0$ and $0 \to (M, \alpha_M) \stackrel{i'} \to (K',\alpha_{K'}) \stackrel{\pi'} \to (L, \alpha_L) \to 0$, are said to be  equivalent if there exists a homomorphism of Hom-Leibniz algebras $\varphi : (K,\alpha_K) \to (K',\alpha_{K'})$ making the following diagram commutative:
\[
\xymatrix{
  0 \ar[r] & (M, \alpha_M) \ar@{=}[d] \ar[r]^{i} & (K,\alpha_K) \ar[d]_{\varphi} \ar[r]^{\pi} & (L, \alpha_L) \ar@{=}[d] \ar[r] & 0  \\
  0 \ar[r] & (M, \alpha_M) \ar[r]^{i'} & (K',\alpha_{K'})  \ar[r]^{\pi'} & (L, \alpha_L) \ar[r] & 0   }
  \]
\end{De}

\begin{Le} \label{ext rota} Let $\left(  C,Id_{C}\right)$ and $\left(  A,\alpha_{A}\right)$ be Hom-Leibniz algebras together with a Hom-action of $\left(  C,Id_{C}\right)$ over $\left(A,\alpha_{A}\right)$.

A sequence of Hom-Leibniz algebras $0\to \left(A,\alpha_{A}\right)  \overset{i}{\to}\left(  B,\alpha_{B}\right)
\overset{\pi}{\to}\left(  C,Id_{C}\right)  \to 0$ is split if and only if it is equivalent to the semi-direct sequence
$0\to \left(  A,\alpha_{A}\right)  \overset{j}{\to}\left(  A\rtimes C,\widetilde{\alpha}\right)  \overset{p}{\to}\left(  C,Id_{C}\right)  \to 0.$
\end{Le}
{\it Proof.}
If $0\to \left(A,\alpha_{A}\right)  \overset{i}{\to}\left(  B,\alpha_{B}\right)
\overset{\pi}{\to}\left(  C,Id_{C}\right)  \to 0$ is split by $s: \left(  C,Id_{C}\right) \to \left(  B,\alpha_{B}\right)$, then the Hom-action of $\left(C,Id_{C}\right)$ over $\left(  A,\alpha_{A}\right)$ is given by
$$c \centerdot  a=\left[  s(c),i(a)\right]; \quad a \centerdot c=\left[i(a),s(c)\right]$$

With this Hom-action of $\left(  C, Id_{C}\right)$ over $\left(  A,\alpha_{A}\right)$ we can construct the following split extension:
\[\xymatrix{ 0 \ar[r] & \left(  A,\alpha_{A}\right) \ar[r]^{j} & \left(  A\rtimes C,\widetilde{\alpha}\right) \ar@<0.5ex>[r]^p &
\left(  C,Id_{C}\right) \ar[r] \ar@<0.5ex>[l]^{\sigma} &0} \]
where $j:A\to  A\rtimes C$, $j(a)=(a,0)$, $p:A\rtimes C \to C,$ $p(a,c)=c$ and $\sigma:C\to A\rtimes C$, $\sigma(c)=(0,c)$.
Moreover the Hom-action of $\left(  C,Id_{C}\right)$ over $\left(  A,\alpha_{A}\right)$ induced by this extension  coincides whit the initial one:
$$c\star a=\left[  \sigma\left(  c\right)  ,j(a)\right]  =\left[  \left(
0,c\right)  \left(  a,0\right)  \right]  =\left(  \left[  0,a\right]
+Id_{C}\left(  c\right)  \centerdot  a+0 \centerdot 0,\left[  c,0\right]  \right)
=(  c \centerdot a, 0) \equiv c \centerdot a$$

Finally, both extensions are  equivalent since the homomorphism of Hom-Leibniz algebras $\varphi:\left(  A\rtimes C,\widetilde{\alpha}\right) \to \left(  B,\alpha_{B}\right)$, $\varphi\left(a,c\right)  =i(a)+s(c)$, makes commutative the following diagram:
\begin{equation} \label{equiv}
 \xymatrix{
0 \ar[r] & \left(  A,\alpha_{A}\right)  \ar[r]^{j} \ar@{=}[d] & \left(  A\rtimes C,\widetilde{\alpha}\right)  \ar@<0.5ex>[r]^{p} \ar@{-->}[d]^{\varphi} & \left(  C,Id_{C}\right) \ar[r] \ar@{=}[d] \ar@<0.5ex>[l]^{\sigma} & 0\\
0 \ar[r] & \left(  A,\alpha_{A}\right)  \ar[r]^{i} & \left(  B,\alpha_{B}\right)  \ar@<0.5ex>[r]^{\pi} & \left(  C,Id_{C}\right) \ar@<0.5ex>[l]^s \ar[r]&  0
}
\end{equation}

For the converse, if both extensions are equivalent, i.e. there exists a homomorphism of Hom-Leibniz algebras $\varphi:\left(  A\rtimes C,\widetilde{\alpha}\right)  \to \left( B,\alpha_{B}\right)$ making commutative diagram (\ref{equiv}), then $s:\left(  C,Id_{C}\right)
\to \left(  B,\alpha_{B}\right)$ given by $s(c)=\varphi\left(0,c\right)$, is a homomorphism that splits the extension. \rdg

\begin{De} Let $\left(  M,\alpha_{M}\right)$ be a  Hom-representation of a Hom-Leibniz algebra  $\left(  L, \alpha_{L}\right)$. A derivation of
$\left(  L, \alpha_{L}\right)$ over $\left(  M,\alpha_{M}\right)$ is a  $\mathbb{K}$-linear map  $d:L\to M$ satisfying:
\begin{enumerate}
\item[a)] $d\left[  l_{1},l_{2}\right]  =\alpha_{L} \left(
l_{1}\right) \centerdot   d\left(  l_{2}\right)  + d\left(  l_{1}\right) \centerdot \alpha
_{L}\left(  l_{2}\right)  $

\item[b)] $d \cdot \alpha_{L}=\alpha_{M} \cdot d$
\end{enumerate}
for all $l_1, l_2 \in L$.
\end{De}

\begin{Ex}\
\begin{enumerate}

\item[a)] The $\mathbb{K}$-linear map $\theta : M \rtimes  L\to M, \theta(m,l)=m,$ is a derivation, where $\left(M,\alpha_{M}\right)$ is a  Hom-representation of $\left(
M \rtimes L,\widetilde{\alpha}\right)$ via $\pi$.

\item[b)] When $(M,\alpha_M)=(L,\alpha_L)$ is considered as a representation following Example \ref{Hom accion Leib} {\it b)}, then a derivation consists of a $\mathbb{K}$-linear map $d:L \to L$ such that $d[l_1,l_2]=[\alpha_L(l_1), d(l_2)]+[d(l_1),\alpha_L(l_2)]$ and $d \cdot \alpha_L = \alpha_L \cdot d$.
\end{enumerate}
\end{Ex}

\begin{Pro}
Let $\left(  M,\alpha_{M}\right)$ be a Hom-representation of a Hom-Leibniz  algebra $\left(  L,\alpha_{L}\right)$.
For every homomorphism of Hom-Leibniz algebras  $f:\left(  X,\alpha_{X}\right)  \to \left(  L,\alpha_{L}\right)$ and every  $f$-derivation
$d:\left(  X,\alpha_{X}\right)  \to \left(  M,\alpha_{M}\right)$ there  exists a unique  homomorphism of Hom-Leibniz algebras  $h:\left(
X,\alpha_{X}\right)  \to \left(  M\rtimes L,\widetilde{\alpha}\right)$, such that the following  diagram is commutative
\[
\xymatrix{
& \left(  X,\alpha_{X} \right) \ar[dr]^{f}  \ar[d]^{h} \ar[ld]^{d} &          \\
 \left(  M,\alpha_{M}\right) \ar[r]^i & \ar@/^0.3pc/[l]^{\theta}     \left(  M\rtimes L,\widetilde{\alpha}\right)  \ar[r]^{\pi} &  \left(  L,\alpha_{L}\right)    }\]

Conversely, every homomorphism of  Hom-Leibniz algebras  $h:\left( X,\alpha_{X}\right)  \to \left(  M\rtimes L,\widetilde{\alpha
}\right)$, determines a homomorphism of Hom-Leibniz algebras   $f=\pi \cdot h: \left(X,\alpha_{X}\right)$  $\to \left(  L,\alpha_{L}\right)$ and any $f$-derivation $d=\theta \cdot h:\left(  X,\alpha_{X}\right)  \to \left(  M,\alpha_{M}\right)$.
\end{Pro}
{\it Proof.} The homomorphism  $h:  X  \to  M\rtimes L, h(x) = \left(  d\left(  x\right)  ,f\left(  x\right)  \right)$ satisfies all the conditions. \rdg

\begin{Co} The set  of all derivations from $\left(  L,\alpha_{L}\right)$ to  $\left(  M,\alpha_{M}\right)$  is in one-to-one correspondence with the set of Hom-Leibniz algebra homomorphisms $h:\left(  L,\alpha_{L}\right)  \to \left(  M\rtimes L,\widetilde{\alpha}\right)$ such that  $\pi \cdot h=Id_L$.
\end{Co}


\section{Functorial properties}

In this section we analyze functorial properties of the universal ($\alpha$)-central extensions of ($\alpha$)-perfect Hom-Leibniz algebras. For detailed motivation, constructions and characterizations we refer to \cite{CIP}.

\begin{De} \label{alfacentral}
A short exact sequence of Hom-Leibniz algebras $(K) : 0 \to (M, \alpha_M) \stackrel{i} \to (K,\alpha_K) \stackrel{\pi} \to (L, \alpha_L) \to 0$ is said to be central if $[M, K] = 0 = [K, M]$. Equivalently,  $M \subseteq Z(K)$.

We say that $(K)$ is $\alpha$-central if $[\alpha_M(M), K] = 0 = [K,\alpha_M(M)]$. Equivalently, $\alpha_M(M) \subseteq Z(K)$.

A central extension $(K) : 0 \to (M, \alpha_M) \stackrel{i} \to (K,\alpha_K) \stackrel{\pi} \to (L, \alpha_L) \to 0$ is said to be universal if for every central extension $(K') : 0 \to (M', \alpha_{M'}) \stackrel{i'} \to (K',\alpha_{K'}) \stackrel{\pi'} \to (L, \alpha_L) \to 0$ there exists a unique homomorphism of Hom-Leibniz algebras  $h : (K,\alpha_K) \to (K',\alpha_{K'})$ such that $\pi'\cdot h = \pi$.

We say that the central extension  $(K) : 0 \to (M, \alpha_M) \stackrel{i} \to (K,\alpha_K) \stackrel{\pi} \to (L, \alpha_L) \to 0$ is universal $\alpha$-central if for every $\alpha$-central extension $(K) : 0 \to (M', \alpha_{M'}) \stackrel{i'} \to (K',\alpha_{K'}) \stackrel{\pi'} \to (L, \alpha_L) \to 0$ there exists a unique homomorphism of Hom-Leibniz algebras  $h : (K,\alpha_K) \to (K',\alpha_{K'})$ such that $\pi'\cdot h = \pi$.
\end{De}

\begin{Rem} \label{rem}
Obviously, every universal $\alpha$-central extension is a universal central extension.
Note that in the case $\alpha_M = Id_M$,  both notions coincide.
\end{Rem}

A perfect ($L = [L,L]$) Hom-Leibniz algebra $(L, \alpha_L)$  admits universal central extension, which is $(\frak{uce}(L), \widetilde{\alpha})$, where  $\frak{uce}(L)=\frac{L \otimes L}{I_L}$ and $I_L$ is the subspace of $L \otimes L$ spanned by the elements of the form $-[x_1,x_2] \otimes \alpha_L(x_3) + [x_1,x_3] \otimes \alpha_L(x_2) + \alpha_L(x_1) \otimes [x_2,x_3], x_1, x_2, x_3 \in L$; every class $x_1 \otimes x_2 + I_L$ is denoted by $\{x_1,x_2\}$, for all $x_1, x_2 \in L$. $\frak{uce}(L)$ is endowed with a structure of Hom-Leibniz algebra with respect to the bracket $[\{x_1,x_2\},\{y_1,y_2\}]=\{[x_1,x_2],[y_1,y_2]\}$  and the endomorphism  $\widetilde{\alpha} : \frak{uce}(L) \to \frak{uce}(L)$  defined by $\widetilde{\alpha}(\{x_1,x_2\}) = \{\alpha_L(x_1), \alpha_L(x_2) \}$.
By construction,  $u_L : (\frak{uce}(L), \widetilde{\alpha}) \to (L,\alpha_L)$, given by $u_L\{x_1,x_2\}=[x_1,x_2]$, gives rise to the universal central extension $0 \to (HL_2^{\alpha}(L), \widetilde{\alpha}_{\mid}) \to (\frak{uce}(L), \widetilde{\alpha}) \stackrel{u_L}\to (L,\alpha_L) \to 0$.

A Hom-Leibniz algebra  $\left(  L, \alpha_{L}\right)$ is said to be $\alpha$-perfect if $L = [\alpha_L(L), \alpha_L(L)]$. Theorem 5.5 in \cite{CIP} shows that a Hom-Leibniz algebra $\left(  L, \alpha_{L}\right)$ is  $\alpha$-perfect if and only if it admits a universal $\alpha$-central extension, which is $(\frak{uce}^{\rm Leib}_{\alpha}(L), \overline{\alpha})$, where $\frak{uce}^{\rm Leib}_{\alpha}(L)= \frac{\alpha_L(L) \otimes \alpha_L(L)}{I_L}$ and $I_L$ is the
vector  subspace spanned by the elements of the form $-[x_1,x_2] \otimes \alpha_L(x_3) + [x_1,x_3] \otimes \alpha_L(x_2) + \alpha_L(x_1) \otimes [x_2,x_3]$, for all $x_1, x_2, x_3 \in L$.
 We denote by $\{\alpha_L(x_1), \alpha_L(x_2)\}$ the equivalence class of $\alpha_L(x_1) \otimes
\alpha_L(x_2) + I_L$. $\frak{uce}^{\rm Leib}_{\alpha}(L)$ is endowed with a structure of Hom-Leibniz algebra with respect to
the bracket $[\{\alpha_L(x_1),\alpha_L(x_2)\}, \{\alpha_L(y_1),\alpha_L(y_2)\}] = \{[\alpha_L(x_1),\alpha_L(x_2)], [\alpha_L(y_1),\alpha_L(y_2)]\}$ and the endomorphism $\overline{\alpha} : \frak{uce}^{\rm Leib}_{\alpha}(L) \to \frak{uce}^{\rm Leib}_{\alpha}(L)$ defined by $\overline{\alpha}(\{\alpha_L(x_1),\alpha_L(x_2)\})$ $= \{\alpha_L^2(x_1),\alpha_L^2(x_2)\}$.
The homomorphism of Hom-Leibniz algebras $U_{\alpha} : \frak{uce}^{\rm Leib}_{\alpha}(L)$ $\to L$
given by $U_{\alpha}(\{\alpha_L(x_1), \alpha_L(x_2)\})= [\alpha_L(x_1), \alpha_L(x_2)]$ gives rise to the universal $\alpha$-central extension $0 \to (Ker (U_{\alpha}), \overline{\alpha}_{\mid}) \to (\frak{uce}_{\alpha}^{\rm Leib}(L), \overline{\alpha})  \stackrel{U_{\alpha}} \to (L, \alpha_L) \to 0$. See \cite{CIP} for details.

\begin{De} A perfect Hom-Leibniz algebra $(L, \alpha_L)$ is said to be centrally closed if its universal central extension is
$$0\to 0\overset{}{\to}\left(  L,\alpha_{L}\right)
\overset{\sim}{\to}\left(  L,\alpha_{L}\right)  \to 0$$
i.e.  $HL_{2}^{\alpha}\left(  L\right)  =0$ and $\left(  \frak{uce}_{\rm Leib} \left( L\right)  ,\widetilde{\alpha}\right)  \cong\left(  L,\alpha_{L}\right)$.

A  Hom-Leibniz algebra $(L, \alpha_L)$ is said to be superperfect if $HL_{1}^{\alpha}\left(  L\right)  = HL_{\substack{2}}^{\alpha
}\left(  L\right)  =0.$
\end{De}

\begin{Co} \label{centralmente cerrada} If $0\to (  Ker(U_{\alpha}), \alpha_{K_{\mid}})   {\to}\left(  K,\alpha_{K}\right)  \stackrel{U_{\alpha}}\to\left(
L,\alpha_{L}\right)  \to 0$ is the universal $\alpha-$central extension of an $\alpha$-perfect Hom-Leibniz algebra $\left(  L,\alpha_{L}\right)$, then $\left(
K,\alpha_{K}\right)$ it is centrally closed.
\end{Co}
{\it Proof.}
By Corollary 4.12 {\it a)} in  \cite{CIP}, $HL_{\substack{1}}^{\alpha}\left(  K\right) = HL_{\substack{2}}^{\alpha}\left(  K\right)=0.$

$HL_{\substack{1}}^{\alpha}\left(  K\right) =0$ if and only if $(K,\alpha_K)$ is perfect. By Theorem 4.11 {\it c)} in \cite{CIP} it admits a universal central extension $0\to (  HL_2^{\alpha}(K), \widetilde{\alpha}_{{\mid}})   {\to}\left( \frak{uce}(K),\widetilde{\alpha} \right)  \stackrel{u}\to\left(
K,\alpha_{K}\right)  \to 0$.  Since $HL_{\substack{2}}^{\alpha}\left(  K\right)=0$, then $u$ is an isomorphism. \rdg

\begin{Le} Let $\pi:\left(  K,\alpha_{K}\right) \twoheadrightarrow \left(  L,\alpha_{L}\right)$ be a  central extension where
$\left(  L,\alpha_{L}\right)$ is a perfect Hom-Leibniz algebra. Then the following statements hold:
\begin{enumerate}
\item[{\it a)}] $K=\left[  K,K\right]  + Ker (\pi)$ and $\overline{\pi}:\left(  \left[
K,K\right],\alpha_{_K{\mid}}\right)  \twoheadrightarrow\left(  L,\alpha_{L}\right)$ is an epimorphism where  $\left(  \left[  K,K\right],\alpha_{\left[K,K\right]}\right)$  is a perfect Hom-Leibniz algebra.

\item[{\it b)}]  $\pi\left(  Z(K)\right)  \subseteq Z(L)$ y $\alpha_L(Z(L)) \subseteq \pi(Z(K))$.
\end{enumerate}
\end{Le}
{\it Proof.}

\noindent {\it a)} It suffices to consider the following commutative diagram:
\[
\xymatrix{
\left(  Ker (\pi) \cap \left[  K,K\right] ,\alpha_{Ker(\pi) \cap\left[  K,K\right] }\right)\  \ar@{>->}[r] \ar@{>->}[d]&
 \left(  \left[  K,K\right]  ,\alpha_{\left[K,K\right]}\right)  \ar@{>>}[r]^{\overline{\pi}} \ar@{>->}[d]&
 \left(\left[  L,L\right],\alpha_{\left[  L,L\right]  }\right) \ar@{=}[d] \\
\left(  Ker (\pi),\alpha_{K\mid}\right)\  \ar@{>>}[d] \ar@{>->}[r] &  \left(  K,\alpha_{K}\right) \ar@{>>}[r]^{\pi} \ar@{>>}[d] & \left(  L,\alpha_{L}\right) \ar@{>>}[d]\\
\ast \ \ar@{>->}[r] &   \left( K/\left[  K,K\right]  ,\overline{\alpha_{K}}\right) \ar@{>>}[r] & \left( L/\left[  L,L\right]  ,\overline{\alpha_{L}}\right)
}
\]

\medskip

\noindent {\it b)}\ Direct checking \rdg

\begin{De} A  Hom-Leibniz algebra $(L, \alpha_L)$ is said to be simply connected if  every central extension $\tau:\left(  F,\alpha_{F}\right) \twoheadrightarrow \left(L,\alpha_{L}\right)$ splits uniquely as the product of  Hom-Leibniz algebras $\left(  F,\alpha_{F}\right)  =\left(  Ker\left(  \tau\right),\alpha_{F\mid}\right)  \times\left(  L,\alpha_{L}\right)$.
\end{De}

\begin{Pro} For a perfect  Hom-Leibniz algebra $\left(  L,\alpha_{L}\right)$,  the following statements are equivalent:
\begin{enumerate}
\item[a)] $\left(  L,\alpha_{L}\right)$ is simply connected.

\item[ b)] $\left(  L,\alpha_{L}\right)$ is centrally closed.

\end{enumerate}

If $u:\left(  L,\alpha_{L}\right)  \twoheadrightarrow\left( M,\alpha_{M}\right)$ is a central extension, then:
  \begin{enumerate}
 \item[c)] Statement {\it a)} (respectively, statement {\it b)}) implies that $u:\left(  L,\alpha_{L}\right)\twoheadrightarrow\left(  M,\alpha_{M}\right)$ is a universal central extension.

\item[d)] If in addition $u:\left(  L,\alpha_{L}\right)  \twoheadrightarrow\left(
M,\alpha_{M}\right)$ is a universal $\alpha$-central extension,  then statements {\it a)} and {\it b)} hold.
\end{enumerate}
\end{Pro}
{\it Proof}.
{\it a)} $\Rightarrow$ {\it b)}   Let $ 0 \to \left(  Ker\left(  u_{\alpha}^{{}}\right)  =HL_{2}^{\alpha}\left(  L\right), \widetilde{\alpha}\right)   \to \left( \frak{uce}_{\alpha}^{{}}\left(  L\right), \widetilde{\alpha}\right) \stackrel{u_{\alpha}}\to \left(  L,\alpha_{L}\right) \to 0$ be the  universal central extension of  $(L, \alpha_L)$, then it is split. Consequently there exists  an isomorphism $\frak{uce}_{\alpha}\left(  L\right)  \cong L$ and $H_{2}^{\alpha}\left( L\right)  =0$.

\medskip

{\it b)} $\Rightarrow$ {\it a)}  The universal central extension of $\left(  L,\alpha_{L}\right)$ is $0\to 0\to \left(  L,\alpha_{L}\right) \overset{\sim}{\to}\left(  L,\alpha_{L}\right)  \to 0$. Consequently every central  extension splits uniquely thanks to the universal property.

\medskip

{\it c)} Let $u:\left(  L,\alpha_{L}\right)  \twoheadrightarrow\left( M,\alpha_{M}\right)$ be a central extension. By Theorem 4.11 {\it b)} in \cite{CIP}, it is  universal if $\left(  L,\alpha_{L}\right)$ is perfect and every central extension of  $\left(  L,\alpha_{L}\right)$ splits.

$\left(  L,\alpha_{L}\right)$ is perfect by hypothesis and by statement {\it a)}, it is simply connected, which means that every central extension splits.

\medskip

{\it d)} If $u:\left(  L,\alpha_{L}\right)  \twoheadrightarrow\left(  M,\alpha_{M}\right)$ is a universal $\alpha-$central extension, then by Theorem 4.1. {\it a)} in \cite{CIP} every central extension $\left(  L,\alpha_{L}\right)$ splits. Consequently  $\left(  L,\alpha_{L}\right)$ is simply connected, equivalently,  it is centrally closed. \rdg

\bigskip

Now we are going to study functorial properties of the universal central extensions.

Consider a homomorphism of perfect Hom-Leibniz algebras $f:\left(  L^{\prime},\alpha_{L^{\prime}}\right)$ $\to
\left( L,\alpha_{L}\right)$. This homomorphism induces a $\mathbb{K}$-linear map $f\otimes f:L^{\prime}\otimes L^{\prime}\to L\otimes L$ given by $\left(  f\otimes f\right)  \left(  x_{1}\otimes x_{2}\right)  = f\left( x_{1}\right)  \otimes f\left(  x_{2}\right)$, that maps the submodule $I_{L'}$ to the submodule $I_L$, hence $f \otimes f$ induces a $\mathbb{K}$-linear map $\frak{uce}(f): \frak{uce}(L^{\prime})\to \frak{uce}(L)$,  given by $\frak{uce}(f)\left\{x_{1},x_{2}\right\}  =\left\{  f(x_{1}),f(x_{2})\right\}$, which is a homomorphism of Hom-Leibniz algebras as well.

Moreover, the following diagram is commutative:
\begin{equation} \label{diagrama uce}
\vcenter{ \xymatrix{
HL_{2}^{\alpha}\left(  L^{\prime}\right) \ar@{>->}[d]   & HL_{2}^{\alpha}\left(  L\right) \ar@{>->}[d] \\
(  \frak{uce}(L^{\prime}),\widetilde{\alpha^{\prime}}) \ar[r]^{\frak{uce}(f)}  \ar@{>>}[d]_{u_{L'}}
& (  \frak{uce}(L),\widetilde{\alpha}) \ar@{>>}[d]^{u_L}\\
(  L^{\prime},\alpha_{L^{\prime}})  \ar[r]^f
& \left(  L,\alpha_{L}\right)
} }\end{equation}

From diagram (\ref{diagrama uce}), the existence of a covariant right exact functor $\frak{uce} :{\sf Hom-Leib^{\rm perf}} \to {\sf Hom-Leib^{\rm perf}}$ between the category of perfect Hom-Leibniz algebras  is derived. Consequently, an automorphism $f$ of $(L,\alpha_L)$ gives rise to an automorphism $\frak{uce}(f)$ of $(\frak{uce}(L), \widetilde{\alpha})$. Commutativity of diagram (\ref{diagrama uce}) implies that $\frak{uce}(f)$ leaves $HL_{2}^{\alpha}(L)$ invariant. So  the  Hom-group homomorphism
$$ \begin{array}{rcl}
{\rm Aut}(L,\alpha_L) & \to & \{ g \in {\rm Aut}(\frak{uce}(L),\widetilde{\alpha}) : f(HL_2^{\alpha}(L))=HL_2^{\alpha}(L) \} \\
f & \mapsto &\frak{uce}(f)
\end{array}
$$ is obtained.

\bigskip

By means of  similar considerations  as the previous ones,  an analogous analysis with respect to the functorial properties of $\alpha$-perfect Hom-Leibniz algebras can be done. Namely, consider a homomorphism of $\alpha-$perfect Hom-Leibniz algebras $f:\left(  L^{\prime},\alpha_{L^{\prime}}\right)  \to \left( L,\alpha_{L}\right)$. Let $I_{L}$ the vector subspace of $\alpha_L(L) \otimes \alpha_L(L)$ spanned by the elements of the form $-[x_1,x_2] \otimes \alpha_L(x_3) + [x_1,x_3] \otimes \alpha_L(x_2) + \alpha_L(x_1) \otimes [x_2,x_3], x_1, x_2,x_2 \in L$, respectively $I_{L^{\prime}}$. $f$ induces a $\mathbb{K}$-linear map  $f\otimes f:\left(  \alpha_{L^{\prime}}\left(  L^{\prime}\right) \otimes\alpha_{L^{\prime}}\left(  L^{\prime}\right) \right.,$ $\left. \alpha_{L^{\prime
}\otimes L^{\prime}}\right)  \to \left(  \alpha_{L}\left(L\right)  \otimes\alpha\left(  _{L}L\right)  ,\alpha_{L\otimes L}\right)$,
given by $\left(  f\otimes f\right)  \left(  \alpha_{L^{\prime}}\left(x_{1}^{\prime}\right) \right.$  $\left. \otimes\alpha_{L^{\prime}}\left(  x_{2}^{\prime}\right)  \right)  =\alpha_L(f\left(  x_{1}^{\prime}\right))  \otimes \alpha_L(f\left(x_{2}^{\prime}\right))$
such that $\left(  f\otimes f\right)  \left(I_{L^{\prime}}\right)  \subseteq I_{L}$.
Consequently, it induces a homomorphism of Hom-Leibniz algebras $\frak{uce}_{\alpha}(f): (\frak{uce}_{\alpha^{\prime}}(L^{\prime}), \overline{\alpha'})\to ( \frak{uce}_{\alpha}(L), \overline{\alpha})$ given by $\frak{uce}_{\alpha}(f)\left\{  \alpha_{L^{\prime}}\left(  x_{1}^{\prime}\right),\alpha_{L^{\prime}}\left(  x_{2}^{\prime}\right)  \right\} =\left\{  \alpha_{L}\left(  f(x'_{1}\right)  ),\alpha_{L}\left(  f(x'_{2})\right)  \right\}$ such that the following diagram is commutative
\begin{equation} \label{diagrama alfa uce}
\vcenter{\xymatrix{
Ker\left(  U_{\alpha^{\prime}}\right) \ \ar@{>->}[d]  & Ker\left(U_{\alpha}\right) \ \ar@{>->}[d] \\
(  \frak{uce}_{\alpha}(L^{\prime}),\overline{\alpha^{\prime}}) \ar[r]^{\frak{uce}_{\alpha}(f)} \ar@{>>}[d]_{U_{\alpha'}} & (  \frak{uce}_{\alpha}(L),\overline{\alpha} )  \ar@{>>}[d]^{U_{\alpha}}\\
(  L^{\prime},\alpha_{L^{\prime}}) \ar[r]^{f}  & (  L,\alpha_{L} )
}}
\end{equation}

From diagram (\ref{diagrama alfa uce}) one derives the existence of a covariant right exact  functor $\frak{uce}_{\alpha} : {\sf Hom-Leib^{\alpha-\rm perf}} \to {\sf Hom-Leib^{\alpha-\rm perf}}$ between the  $\alpha$-perfect Hom-Leibniz algebras category. Consequently, an automorphism $f$ of $(L,\alpha_L)$ gives rise to an automorphism $\frak{uce}_{\alpha}(f)$ of $(\frak{uce}_{\alpha}(L), \overline{\alpha})$. Commutativity of diagram (\ref{diagrama alfa uce}) implies that $\frak{uce}_{\alpha}(f)$ leaves $Ker(U_{\alpha})$ invariant. So the homomorphism of Hom-groups
$$ \begin{array}{rcl}
{\rm Aut}(L,\alpha_L) & \to & \{ g \in {\rm Aut}(\frak{uce}_{\alpha}(L),\overline{\alpha}) : f(Ker(U_{\alpha})=Ker(U_{\alpha}) \} \\
f & \mapsto &\frak{uce}_{\alpha}(f)
\end{array}
$$
is obtained.

\bigskip

Now we consider a derivation $d$ of the   $\alpha$-perfect Hom-Leibniz algebra $\left(  L,\alpha_{L}\right)$. The linear map $\varphi:\alpha_L(L)\otimes \alpha_L(L) \to  \alpha_L(L)\otimes \alpha_L(L)$ given by $\varphi\left(  \alpha_L(x_{1})\otimes \alpha_L(x_{2})\right)  =d(\alpha_L(x_{1}))\otimes\alpha^2 _{L}\left(  x_{2}\right)  +\alpha^2_{L}\left(  x_{1}\right)  \otimes d\left(\alpha_L(x_{2})\right)$, leaves invariant the  vector subspace $I_{L}$ of $\alpha_L(L) \otimes \alpha_L(L)$ spanned by the  elements of the form $ -[x_1,x_2] \otimes \alpha_L(x_3) + [x_1,x_3] \otimes \alpha_L(x_2) + \alpha_L(x_1) \otimes [x_2,x_3], x_1, x_2,x_2 \in L$. Hence it induces a linear map $\frak{uce}_{\alpha}(d):\left(  \frak{uce}_{\alpha}\left(  L\right),\overline{\alpha}\right)  \to \left( \frak{uce}_{\alpha}\left(  L\right),\overline{\alpha}\right)$, given by
$\frak{uce}_{\alpha}(d)( \left\{  \alpha_L(x_{1}),\alpha_L(x_{2}) \right\}) = \left\{  d(\alpha_L(x_{1})),\alpha^2_{L}\left(
x_{2}\right)  \right\}  +\left\{  \alpha^2_{L}\left(  x_{1}\right)  ,d\left(
\alpha_L(x_{2})\right)  \right\}$, that makes commutative the following diagram:
\begin{equation} \label{diagrama derivacion}
\vcenter{ \xymatrix{
\left(  \frak{uce}_{\alpha}\left(  L\right),\overline{\alpha} \right)  \ar[r]^{uce_{\alpha}(d)} \ar@{>>}[d]_{U_{\alpha}} & \left(  \frak{uce}_{\alpha}\left(  L\right),\overline{\alpha} \right)  \ar@{>>}[d]^{U_{\alpha}}  \\
\left(  L,\alpha_{L}\right)  \ar[r]^d & \left(  L,\alpha_{L}\right)
} }
\end{equation}
Consequently, a derivation $d$ of $(L,\alpha_L)$ gives rise to a derivation $\frak{uce}_{\alpha}(d)$ of $(\frak{uce}_{\alpha}(L), \overline{\alpha})$. The commutativity of diagram (\ref{diagrama derivacion}) implies that $\frak{uce}_{\alpha}(d)$ maps $Ker(U_{\alpha})$ on itself.

Hence, it is obtained  the homomorphism of Hom-$\mathbb{K}$-vector spaces
$$ \begin{array}{rcl}
\frak{uce}_{\alpha} : {\rm Der}(L,\alpha_L) & \to & \{ \delta \in {\rm Der}(\frak{uce}_{\alpha}(L),\overline{\alpha}) : \delta (Ker(U_{\alpha}) \subseteq Ker(U_{\alpha}) \} \\
d & \mapsto &\frak{uce}_{\alpha}(d)
\end{array}
$$
whose kernel belongs to the subalgebra of derivations of $(L,\alpha_L)$ such that vanish on $[\alpha_L(L),\alpha_L(L)]$.

The functorial properties of  $\frak{uce}_{\alpha}(-)$ relative to the  derivations are described by the following result.

\begin{Le} \label{uce derivacion}
Let  $f:(L^{\prime},\alpha_{L^{\prime}})\to \left(  L,\alpha_{L}\right)$ be a homomorphism of  $\alpha$-perfect Hom-Leibniz algebras.
Consider $d\in Der(L)$ and  $d^{\prime}\in Der(L^{\prime})$ such that
$f \cdot d^{\prime}=d \cdot f$. then $\frak{uce}_{\alpha}(f) \cdot \frak{uce}_{\alpha}(d^{\prime})= \frak{uce}_{\alpha}(d) \cdot \frak{uce}_{\alpha}(f).$
\end{Le}
{\it Proof.} Routine checking. \rdg

\section{Lifting automorphisms and derivations}

In this section we analyze under what conditions an automorphism or a derivation can be lifted to an $\alpha$-cover. We restrict the study to $\alpha$-covers since we must compose central extensions in the  constructions. This fact does not allow to obtain more general results, mainly due to Lemma 4.10 in \cite{CIP}.

\begin{De} A central extension of Hom-Leibniz algebras $f:\left(  L^{\prime},\alpha_{L^{\prime}}\right)$
 $\twoheadrightarrow\left(  L,\alpha_{L}\right)$, where  $\left(  L^{\prime},\alpha_{L^{\prime}}\right)$ is an $\alpha$-perfect Hom-Leibniz  algebra, is said to be an $\alpha$-cover.
\end{De}

\begin{Le} \label{alfa perfecta sobre} If $f:\left(  L^{\prime},\alpha_{L^{\prime}}\right) \to \left(  L,\alpha_{L}\right)$ is a surjective homomorphism of Hom-Leibniz algebras and  $\left(  L^{\prime},\alpha_{L^{\prime}}\right)$ is $\alpha-$perfect,
then $\left(  L,\alpha_{L}\right)$ is $\alpha$-perfect as well.
\end{Le}

Let $f:\left(  L^{\prime},\alpha_{L^{\prime}}\right)\twoheadrightarrow \left(  L,\alpha_{L}\right)$ be an $\alpha$-cover.
 Thanks to Lema  \ref{alfa perfecta sobre} $\left(  L,\alpha_{L}\right)$  is an $\alpha$-perfect  Hom-Leibniz algebra as well.  By Theorem 5.5 in \cite{CIP}, everyone admits  universal $\alpha$-central extension. Having in mind the functorial properties given in diagram (\ref{diagrama alfa uce}), we can construct the following diagram:

\[{\xymatrix{
Ker\left(  U_{\alpha^{\prime}}\right) \ \ar@{>->}[d]  & Ker\left(U_{\alpha}\right) \ \ar@{>->}[d] \\
(  \frak{uce}_{\alpha}(L^{\prime}),\overline{\alpha^{\prime}}) \ar[r]^{\frak{uce}_{\alpha}(f)} \ar@{>>}[d]_{U_{\alpha'}} & (  \frak{uce}_{\alpha}(L),\overline{\alpha} )  \ar@{>>}[d]^{U_{\alpha}}\\
(  L^{\prime},\alpha_{L^{\prime}}) \ar[r]^{f}  & (  L,\alpha_{L} )
}}\]

Since $U_{\alpha^{\prime}}:\left(  \frak{uce}_{\alpha^{\prime}}\left(  L^{\prime}\right),\overline{\alpha^{\prime}}\right)
\twoheadrightarrow \left(  L^{\prime},\alpha_{L^{\prime}}\right)$ is a universal $\alpha$-central extension, then by Remark \ref{rem}, it  is a universal  central extension as well.
Since  $f:\left(  L^{\prime},\alpha_{L^{\prime}}\right)  \twoheadrightarrow \left(  L,\alpha_{L}\right)$ is a central extension and   $U_{\alpha^{\prime}}: \left(  \frak{uce}_{\alpha^{\prime}}\left( L^{\prime}\right),\overline{\alpha^{\prime}}\right)  \twoheadrightarrow \left(  L',\alpha_{L'}\right)$ is a universal central extension, then by Proposition 4.15 in \cite{CIP} the extension $f \cdot U_{\alpha^{\prime}}:\left(  \frak{uce}_{\alpha^{\prime}}\left(  L^{\prime}\right)  ,\widetilde{\alpha}^{\prime}\right)  \to\left(  L,\alpha_{L}\right)$ is $\alpha$-central which is universal in the sense of Definition 4.13 in \cite{CIP}.

On the other hand, since $U_{\alpha} : \left(  \frak{uce}_{\alpha}\left(  L\right),\overline{\alpha} \right) \twoheadrightarrow (L,\alpha_L)$ is a universal $\alpha$-central extension, then there exists a unique homomorphism $\varphi:\left(  \frak{uce}_{\alpha}\left(  L\right), \overline{\alpha} \right)  \twoheadrightarrow \left(  \frak{uce}_{\alpha^{\prime}}\left(L^{\prime}\right),\overline{\alpha^{\prime}}\right)$ such that $f \cdot U_{\alpha^{\prime}} \cdot \varphi=U_{\alpha}$.

Moreover $\varphi \cdot \frak{uce}_{\alpha}(f) = Id$ since the following  diagram is commutative
\[\xymatrix{
0 \ar[r] & \left(  Ker(f \cdot U_{\alpha^{\prime}}),\overline{\alpha'_{\mid}}\right)  \ar[r]& \left( \frak{uce}_{\alpha^{\prime}}\left(  L^{\prime}\right),\overline{\alpha'}\right)  \ar[r]^{\quad f \cdot U_{\alpha^{\prime}}} \ar@<2ex>[d]_{\varphi \cdot \frak{uce}_{\alpha}(f)\quad} \ar[d]^{\quad Id}& \left(  L,\alpha_{L}\right) \ar[r] \ar@{=}[d]& 0 \\
0 \ar[r] & \left(  Ker(f \cdot U_{\alpha^{\prime}}),\overline{\alpha'_{\mid}}\right) \ar[r] & \left(  \frak{uce}_{\alpha^{\prime}}\left(  L^{\prime}\right),\overline{\alpha^{\prime}}\right)  \ar[r]^{\quad f \cdot U_{\alpha^{\prime}}} & \left(  L,\alpha_{L}\right) \ar[r] & 0
}\]
and $f \cdot U_{\alpha'}$ is an $\alpha$-central extension which is universal in the sense of Definition 4.13 in \cite{CIP}.

Conversely,  $\frak{uce}_{\alpha}(f) \cdot \varphi = Id$ since the following diagram is commutative
\[\xymatrix{
0 \ar[r] & \left(  Ker( U_{\alpha}),\overline{\alpha_{\mid}}\right)  \ar[r]& \left( \frak{uce}_{\alpha}\left(  L\right),\overline{\alpha}\right)  \ar[r]^{\quad U_{\alpha}} \ar@<2ex>[d]_{\frak{uce}_{\alpha}(f) \cdot \varphi\quad} \ar[d]^{\quad Id}& \left(  L,\alpha_{L}\right) \ar[r] \ar@{=}[d]& 0 \\
0 \ar[r] & \left(  Ker(U_{\alpha}),\overline{\alpha_{\mid}}\right) \ar[r] & \left(  \frak{uce}_{\alpha}\left(  L\right),\overline{\alpha}\right)  \ar[r]^{\quad  U_{\alpha}} & \left(  L,\alpha_{L}\right) \ar[r] & 0
}\]
whose horizontal rows are central extensions and $\left( \frak{uce}_{\alpha}\left(  L\right),\overline{\alpha}\right)$ is $\alpha$-perfect, then Lemma 5.4 in \cite{CIP} guarantees the uniqueness of the vertical  homomorphism.

Consequently $\frak{uce}_{\alpha}(f)$ is an isomorphism and from now on we will use  the notation $\frak{uce}_{\alpha}(f)^{-1}$ instead of $\varphi$.

On the other hand, $U_{\alpha^{\prime}} \cdot \frak{uce}_{\alpha}(f)^{-1}:\left(  \frak{uce}_{\alpha}\left(  L\right)  ,\overline{\alpha} \right)$  $\twoheadrightarrow \left(  L^{\prime},\alpha_{L^{\prime}}\right)$ is an $\alpha-$cover. In the sequel, we will denote its kernel by
$$C:=Ker(U_{\alpha^{\prime}}\cdot \frak{uce}_{\alpha}(f)^{-1})= \frak{uce}_{\alpha}(f)\left(  Ker\left(  U_{\alpha^{\prime}}\right)  \right)$$

\begin{Th} \label{levantamiento automorfismo} Let $f:\left(  L^{\prime},\alpha_{L^{\prime}}\right) \twoheadrightarrow \left(  L,\alpha_{L}\right)$ be an $\alpha-$cover.

For any $h\in Aut\left(  L,\alpha_{L}\right)$, there exists a unique $\theta_{h}\in Aut\left(  L^{\prime},\alpha_{L^{\prime}}\right)$ such that the following diagram is commutative:
\begin{equation} \label{automorfismo}
\vcenter{ \xymatrix{
\left(  L^{\prime},\alpha_{L^{\prime}}\right)  \ar@{>>}[r]^{f} \ar[d]_{\theta_h}& \left(  L,\alpha_{L}\right) \ar[d]^h \\
\left(  L^{\prime},\alpha_{L^{\prime}}\right)  \ar@{>>}[r]^{f} & \left(  L,\alpha_{L}\right)
}}
\end{equation}
if and only if the automorphism $\frak{uce}_{\alpha}(h)$ of $(\frak{uce}_{\alpha}\left( L \right), \overline{\alpha})$ satisfies $\frak{uce}_{\alpha}(h)\left(  C\right)  =C$. In this case,  it is uniquely determined by diagram (\ref{automorfismo}) and $\theta_{h}\left(  Ker(f)\right)  =Ker(f)$.

Moreover, the map
$$\begin{array}{rcl}\Theta:\left\{  h\in Aut\left(  L,\alpha_{L}\right)  : \frak{uce}_{\alpha}(h)\left( C\right)  =C\right\}  &\to & \left\{  g\in Aut\left(  L^{\prime},\alpha_{L^{\prime}}\right)  :g\left(  Ker(f)\right)  =Ker(f)\right\}\\
 h &\mapsto & \theta_{h} \end{array}$$
is a group isomorphism.
\end{Th}
{\it Proof.} Let $h\in Aut\left(  L,\alpha_{L}\right)$ be and  assume that there exists a
$\theta_{h}\in Aut\left(  L^{\prime},\alpha_{L^{\prime}}\right)$ such that diagram (\ref{automorfismo}) is commutative.

Then $h \cdot f:\left(  L^{\prime},\alpha_{L^{\prime}}\right)  \to\left(  L,\alpha_{L}\right)$ is an $\alpha-$cover, hence $\theta_h$ is a homomorphism from the $\alpha$-cover $h \cdot f$ to the $\alpha$-cover $f$ which is unique by Remark 5.3 {\it b)} and Lemma 4.7 in \cite{CIP}.

By application of the functor $\frak{uce}_{\alpha}(-)$ to diagram (\ref{automorfismo}), one obtains the following commutative diagram:
\[\xymatrix{
\left(  \frak{uce}_{\alpha}\left(  L^{\prime}\right),\overline{\alpha_L{^{\prime}}}\right) \ar@{>>}[r]^{\frak{uce}_{\alpha}(f)} \ar@{>->>}[d]_{\frak{uce}_{\alpha}(\theta_{h})}& \left(  \frak{uce}_{\alpha}(L),\overline{\alpha_L} \right) \ar@{>->>}[d]^{\frak{uce}_{\alpha}(h)}\\
\left(  \frak{uce}_{\alpha}\left(  L^{\prime}\right),\overline{\alpha_{L^{\prime}}}\right)  \ar@{>>}[r]^{\frak{uce}_{\alpha}(f)}
 & \left(  \frak{uce}_{\alpha}(L),\overline{\alpha_L} \right)
}\]

Hence $\frak{uce}_{\alpha}(h)\left(  C\right)  = \frak{uce}_{\alpha}(h)\cdot \frak{uce}_{\alpha}(f)\left(  Ker\left(  U_{\alpha'}\right)  \right) =\frak{uce}_{\alpha}(f) \cdot \frak{uce}_{\alpha}(\theta_{h})\left(  Ker\left(  U_{\alpha'}\right) \right) = \frak{uce}_{\alpha}(f) \left(  Ker\left(  U_{\alpha'}\right) \right) = C$

Conversely, from diagram (\ref{diagrama alfa uce}), we have that $U_{\alpha}=f \cdot U_{\alpha^{\prime}} \cdot \frak{uce}_{\alpha}(f)^{-1}$, hence we obtain the following diagram:
\[\xymatrix{
C \ar@{>->}[r] \ar[d] & \left(  \frak{uce}_{\alpha}(L),\overline{\alpha} \right) \ar@{>>}[rr]^{U_{\alpha^{\prime}}\cdot \frak{uce}_{\alpha}(f)^{-1}} \ar[d]^{\frak{uce}_{\alpha}(h)}& \ & \left(L^{\prime},\alpha_{L^{\prime}}\right)  \ar@{>>}[r]^{f} \ar@{-->}[d]^{\theta_h} & \left(  L,\alpha_{L}\right) \ar[d]^h\\
C \ar@{>->}[r] & \left(  \frak{uce}_{\alpha}(L),\overline{\alpha} \right) \ar@{>>}[rr]^{U_{\alpha^{\prime}}\cdot \frak{uce}_{\alpha}(f)^{-1}} & \ & \left(L^{\prime},\alpha_{L^{\prime}}\right)  \ar@{>>}[r]^{f}  & \left(  L,\alpha_{L}\right)
}\]
 If $\frak{uce}_{\alpha}(h)\left(  C\right)  =C$, then $U_{\alpha^{\prime}}\cdot \frak{uce}_{\alpha}(f)^{-1}\cdot \frak{uce}_{\alpha}(h)\left(  C\right)  = U_{\alpha^{\prime}}\cdot \frak{uce}_{\alpha}(f)^{-1}(C)=0$, then there  exists a unique
$\theta_{h}:\left(  L^{\prime},\alpha_{L^{\prime}}\right)  \to \left(  L^{\prime},\alpha_{L^{\prime}}\right)$ such that $\theta_{h}
\cdot U_{\alpha^{\prime}} \cdot \frak{uce}_{\alpha}(f)^{-1}=U_{\alpha^{\prime}} \cdot \frak{uce}_{\alpha}(f)^{-1} \cdot \frak{uce}_{\alpha}(h)$.

On the other hand, $h\cdot f \cdot U_{\alpha^{\prime}} \cdot \frak{uce}_{\alpha}(f)^{-1}=f \cdot U_{\alpha^{\prime}} \cdot \frak{uce}_{\alpha}(f)^{-1} \cdot \frak{uce}_{\alpha}(h)= f \cdot U_{\alpha'} \cdot \frak{uce}_{\alpha}(\theta_h) \cdot \frak{uce}_{\alpha}(f)^{-1} = f \cdot \theta_{h} \cdot U_{\alpha^{\prime}} \cdot \frak{uce}_{\alpha}(f)^{-1}$, then $h\cdot f=f\cdot \theta_{h}.$

In conclusion, $\theta_{h}$ is uniquely determined by diagram (\ref{automorfismo}) and moreover $\theta_{h}( Ker(f))$
 $=Ker(f)$.

By the previous arguments, it is easy to check that $\Theta$ is a well-defined map, it is a monomorphism thanks to the uniqueness of $\theta_{h}$ and it is an epimorphism, since every $g\in Aut\left(  L',\alpha_{L'}\right)$ with $g\left(  Ker(f)\right)
=Ker(f)$, induces  a unique homomorphism $h:\left(  L,\alpha_{L}\right)  \to \left(  L,\alpha_{L}\right)$ such that $h \cdot f=f \cdot g$. Then $g=\theta_{h}$ and $\frak{uce}_{\alpha}\left(h\right)  \left(  C\right)  =C$. \rdg

\begin{Co} If $\left(  L,\alpha_{L}\right)$ is an $\alpha$-perfect Hom-Leibniz algebra, then the map
$$\begin{array}{rcl} Aut(L,\alpha_{L}) & \to &\left\{  g\in Aut(\frak{uce}_{\alpha}(L),\overline{\alpha}):g(Ker(U_{\alpha}))=Ker(U_{\alpha})\right\}\\
h &\mapsto & \frak{uce}_{\alpha}(h)\end{array}$$
is a group isomorphism.
\end{Co}
{\it Proof}. By application of Theorem \ref{levantamiento automorfismo} to the $\alpha-$cover $U_{\alpha}:\left(  \frak{uce}_{\alpha
}\left(  L\right),\overline{\alpha}\right)  \twoheadrightarrow\left(L,\alpha_{L}\right)$, it is enough to have in mind that under these conditions $C=0$  and $\frak{uce}_{\alpha}(f)(0)=0$. \rdg
\bigskip

Now we analyze under what conditions a derivation of an $\alpha$-perfect Hom-Leibniz algebra can be lifted to an $\alpha$-cover.

\begin{Th} Let $f:\left(  L^{\prime},\alpha_{L^{\prime}}\right)  \twoheadrightarrow\left(  L,\alpha_{L}\right)$ be an $\alpha-$cover.
Denote by $C=\frak{uce}_{\alpha}(f)\left(  Ker\left(  U_{\alpha^{\prime}}\right)  \right)  \subseteq Ker(U_{\alpha})$. Then the following statements hold:
\begin{enumerate}
\item[a)] For any $d\in Der\left(  L,\alpha_{L}\right)$ there exists a
$\delta_{d}\in Der\left(  L^{\prime},\alpha_{L^{\prime}}\right)$ such that the following diagram is commutative
 \begin{equation} \label{derivacion}
\vcenter{ \xymatrix{
\left(  L^{\prime},\alpha_{L^{\prime}}\right)  \ar@{>>}[r]^{f} \ar[d]_{\delta_d}& \left(  L,\alpha_{L}\right) \ar[d]^d \\
\left(  L^{\prime},\alpha_{L^{\prime}}\right)  \ar@{>>}[r]^{f} & \left(  L,\alpha_{L}\right)
}}
\end{equation}
if and only if the derivation $\frak{uce}_{\alpha}(d)$ of $(\frak{uce}_{\alpha}\left(L\right),\overline{\alpha_L})$ satisfies $\frak{uce}_{\alpha}(d)\left(  C\right)  \subseteq C$.

In this case, $\delta_{d}$ is uniquely determined by (\ref{derivacion}) and $\delta_{d}\left(  Ker(f)\right)$  $\subseteq  Ker(f)$.

\item[b)] The map
$$ \begin{array}{rcl} \Delta:\left\{  d\in Der\left(  L,\alpha_{L}\right)
:\frak{uce}_{\alpha}(d)\left(  C\right)  \subseteq C\right\}  & \to &\left\{\rho\in Der\left(  L^{\prime}, \alpha_{L^{\prime}}\right)  :\rho\left(Ker(f)\right)  \subseteq Ker(f)\right\}\\
 d & \mapsto & \delta_{d}\end{array}$$
is an isomorphism of Hom-vector spaces.

\item[c)] For the $\alpha$-cover $U_{\alpha} : (\frak{uce}_{\alpha}(L),\overline{\alpha_L}) \twoheadrightarrow (L, \alpha_L)$, the map
$$\frak{uce}_{\alpha} : Der(L,\alpha_L) \to \{\delta \in Der (\frak{uce}_{\alpha}(L),\overline{\alpha_L}) : \delta(Ker(U_{\alpha})) \subseteq Ker(U_{\alpha}) \}$$
is an isomorphism of Hom-vector spaces.
\end{enumerate}
\end{Th}
{\it Proof}. {\it a)} Let $d\in Der\left(  L, \alpha_L \right)$ be and assume the  existence of a $\delta_{d}\in
Der\left(  L', \alpha_{L'} \right)$ such that  diagram(\ref{derivacion}) is commutative. Then, by Lemma \ref{uce derivacion}, we obtain the following diagram commutative:
\[\xymatrix{
\left(  \frak{uce}_{\alpha}\left(  L^{\prime}\right),\overline{\alpha_{L^{\prime}}}\right)  \ar[r]^{\frak{uce}_{\alpha}(f)} \ar[d]_{\frak{uce}_{\alpha}(\delta_d)} & \left(\frak{uce}_{\alpha}(L),\overline{\alpha_L} \right) \ar[d]^{\frak{uce}_{\alpha}(d)}\\
\left(  \frak{uce}_{\alpha}\left(  L^{\prime}\right), \overline{\alpha_{L^{\prime}}} \right)  \ar[r]^{\frak{uce}_{\alpha}(f)} & \left(
\frak{uce}_{\alpha}(L),\overline{\alpha_L} \right)
}\]
Hence, having in mind the properties derived from diagram (\ref{diagrama derivacion}), we obtain:

\noindent $\frak{uce}_{\alpha}(d)\left(  C\right)  = \frak{uce}_{\alpha}(d) \cdot \frak{uce}_{\alpha}(f) \left(  Ker(U_{\alpha^{\prime}})\right)  =\frak{uce}_{\alpha}(f) \cdot \frak{uce}_{\alpha}(\delta_d)\left(  Ker(U_{\alpha^{\prime}})\right) \subseteq \frak{uce}_{\alpha}(f)\left(  Ker(U_{\alpha^{\prime}})\right) = C$.

Conversely, from diagram (\ref{diagrama alfa uce}) we have that $U_{\alpha}=f \cdot U_{\alpha^{\prime}} \cdot \frak{uce}_{\alpha}(f)^{-1}$ and consider the following  diagram:
\[\xymatrix{
C \ar@{>->}[r] \ar[d] & \left(  \frak{uce}_{\alpha}(L),\overline{\alpha_L} \right) \ar@{>>}[rr]^{U_{\alpha^{\prime}}\cdot \frak{uce}_{\alpha}(f)^{-1}} \ar[d]^{\frak{uce}_{\alpha}(d)}& \ & \left(L^{\prime},\alpha_{L^{\prime}}\right)  \ar@{>>}[r]^{f} \ar@{-->}[d]^{\delta_d} & \left(  L,\alpha_{L}\right) \ar[d]^d\\
C \ar@{>->}[r] & \left(  \frak{uce}_{\alpha}(L),\overline{\alpha_L} \right) \ar@{>>}[rr]^{U_{\alpha^{\prime}}\cdot \frak{uce}_{\alpha}(f)^{-1}} & \ & \left(L^{\prime},\alpha_{L^{\prime}}\right)  \ar@{>>}[r]^{f}  & \left(  L,\alpha_{L}\right)
}\]
Since $\frak{uce}_{\alpha}(d)\left(  C\right)  \subseteq C$, then
$U_{\alpha^{\prime}} \cdot \frak{uce}_{\alpha}(f)^{-1} \cdot \frak{uce}_{\alpha}(d)\left(  C\right)
\subseteq U_{\alpha^{\prime}} \cdot \frak{uce}_{\alpha}(f)^{-1}(C)= U_{\alpha'} (Ker(U_{\alpha'})) = 0$.
Hence there exists a unique $\mathbb{K}$-linear map  $\delta_{d}:\left(  L^{\prime},\alpha_{L^{\prime}}\right)  \to \left(  L^{\prime},\alpha_{L^{\prime}}\right)$ such that $\delta_{d} \cdot U_{\alpha^{\prime}} \cdot \frak{uce}_{\alpha}(f)^{-1}=U_{\alpha^{\prime}} \cdot \frak{uce}_{\alpha}(f)^{-1} \cdot \frak{uce}_{\alpha}(d).$

 On the other hand $d \cdot f \cdot U_{\alpha^{\prime}} \cdot \frak{uce}_{\alpha}(f)^{-1} = d \cdot U_{\alpha} \cdot \frak{uce}_{\alpha}(f) \cdot \frak{uce}_{\alpha}(f)^{-1} = U_{\alpha} \cdot \frak{uce}_{\alpha}(d) =f \cdot U_{\alpha^{\prime}} \cdot \frak{uce}_{\alpha}(f)^{-1} \cdot \frak{uce}_{\alpha}(d)$, since $\delta_d \cdot U_{\alpha'} \cdot \frak{uce}_{\alpha}(f)^{-1} = U_{\alpha'} \cdot \frak{uce}_{\alpha}(f)^{-1} \cdot \frak{uce}_{\alpha}(d)$, then $d \cdot f=f \cdot \delta_{d}$.

Finally, a direct checking shows that $\delta_{d}$ is a derivation of $L^{\prime}$, which is uniquely determined by diagram (\ref{derivacion}) and
$\delta_{d}\left(  Ker(f)\right)  \subseteq Ker(f)$.

\bigskip

{\it b)} The map $\Delta$ is a homomorphism of Hom-vector spaces by construction, which is injective by the uniqueness of $\delta_{d}$, and surjective, since for every $\rho\in Der\left(  L^{\prime},\alpha_{L^{\prime}}\right)$ such that  $\rho\left(  Ker(f)\right)  \subseteq Ker(f)$ there exists the following diagram commutative:
\[ \xymatrix{
Ker(f)\ \ar@{>->}[r] \ar[d] & \left(  L^{\prime},\alpha_{L^{\prime}}\right)  \ar@{>>}[r]^f \ar[d]^{\rho} & \left(  L,\alpha_{L}\right) \ar@{-->}[d]^d\\
Ker(f)\ \ar@{>->}[r] & \left(  L^{\prime},\alpha_{L^{\prime}}\right)  \ar@{>>}[r]^f & \left(  L,\alpha_{L}\right)
} \]
where $d:\left(  L,\alpha_{L}\right)  \to \left(  L,\alpha_{L}\right)$ is a derivation  satisfying
$\frak{uce}_{\alpha}\left(  d\right)  \left(  C\right)  = \frak{uce}_{\alpha}\left(
d\right)  \cdot \frak{uce}_{\alpha}(f)\left(  Ker(U_{\alpha^{\prime}})\right)
=\frak{uce}_{\alpha}\left(  f\right)  \cdot \frak{uce}_{\alpha}(\rho)\left(  Ker(U_{\alpha^{\prime}})\right)  \subseteq \frak{uce}_{\alpha}(f)\left(  Ker(U_{\alpha^{\prime}})\right)  =C.$

Finally,   the uniqueness of $\delta_d$ implies that $\Delta\left(  d\right)  = \delta_{d}= \rho$.

\bigskip

{\it c)} It is enough to write the statement   {\it b)} for the $\alpha-$cover
$U_{\alpha}:\left(  \frak{uce}_{\alpha}(L),\overline{\alpha_{L}}\right)\twoheadrightarrow\left(  L,\alpha_{L}\right)$. Now $C= \frak{uce}_{\alpha}(U_{\alpha})\left(  Ker\left(  U_{\alpha}\right)  \right)  =0$, and $\Delta$ is the map $\frak{uce}_{\alpha}$  derived from diagram (\ref{diagrama derivacion}). \rdg

\section{Universal $\alpha$-central extension of the semi-direct product}

Consider a split extension of  $\alpha$-perfect Hom-Leibniz algebras
\[\xymatrix{
0 \ar[r]& (M,\alpha_M) \ar[r]^t& (G,\alpha_G) \ar@<0.5ex>[r]^p &(Q,Id_Q) \ar[r] \ar@<0.5ex>[l]^s &0} \]
where, by Lemma \ref{ext rota},  $(G,\alpha_G) \cong (M,\alpha_M) \rtimes (Q,Id_Q)$, whose Hom-action of $(Q,Id_Q)$ on $(M,\alpha_M)$ is given by $q \centerdot m = [s(q), t(m)]$ and $m \centerdot q = [t(m),s(q)], q \in Q, m \in M$. Moreover we will assume, when it is needed, that the previous action  is symmetric, i.e. $q \centerdot m + m \centerdot q =0, q \in Q, m \in M$.

An example of the above situation is when $(M,\alpha_M)$ is an $\alpha$-perfect Hom-Leibniz algebra, $Q$ is a perfect Leibniz algebra considered as the Hom-Leibniz algebra $(Q, Id_Q)$ and $(G,\alpha_G)=(M,\alpha_M) \times (Q, Id_Q)=(M \times Q, \alpha_M \times Id_Q)$.

Applying the functorial properties of $\frak{uce}_{\alpha}(-)$ given in diagram (\ref{diagrama alfa uce}) and having in mind  that $(Q,Id_Q)$ is perfect is equivalent to $Q$ is perfect, we have the following commutative diagram:
\[ \xymatrix{
& Ker(U_{\alpha}^M) \ar@{>->}[d]& Ker(U_{\alpha}^G) \ar@{>->}[d]& HL_2(Q) \ar@{>->}[d]& \\
& (\frak{uce}_{\alpha}(M),\overline{\alpha_M}) \ar@{>>}[d]^{U_{\alpha}^M}\ar[r]^{\tau}& (\frak{uce}_{\alpha}(G),\overline{\alpha_G}) \ar@{>>}[d]^{U_{\alpha}^G} \ar@<0.5ex>[r]^{\pi} &(\frak{uce}(Q),Id_{\frak{uce}(Q)}) \ar@{>>}[d]^{u_Q} \ar@<0.5ex>[l]^{\sigma} &\\
0 \ar[r]& (M,\alpha_M) \ar[r]^t& (G,\alpha_G) \ar@<0.5ex>[r]^p &(Q,Id_Q) \ar[r] \ar@<0.5ex>[l]^s & 0
} \]
Here $\tau = \frak{uce}_{\alpha}(t), \pi = \frak{uce}_{\alpha}(p), \sigma = \frak{uce}_{\alpha}(s)$.

The sequence \[\xymatrix{
& (\frak{uce}_{\alpha}(M),\overline{\alpha_M}) \ar[r]^{\tau}& (\frak{uce}_{\alpha}(G),\overline{\alpha_G})  \ar@<0.5ex>[r]^{\pi} &(\frak{uce}(Q),Id_{\frak{uce}(Q)})  \ar@<0.5ex>[l]^{\sigma} &} \]
is split, since $p \cdot s = Id_Q$, then $\frak{uce}_{\alpha}(p) \cdot \frak{uce}_{\alpha}(s) = \frak{uce}_{\alpha}(Id_Q)$, i.e. $\pi \cdot \sigma = Id_{\frak{uce}(Q)}$. Obviously  $\pi$ is surjective and there exists a Hom-action of $(\frak{uce}(Q),Id_{\frak{uce}(Q)})$ on $(Ker(\pi),\overline{\alpha_{G}}_{\mid})$ induced by the section $\sigma$, which is given by:

\noindent  $\lambda:\frak{uce}(Q)\otimes Ker(\pi) \to Ker(\pi)$,

 $\lambda\left(  \left\{ q_{1}, q_{2}  \right\}  \otimes \left\{  \alpha_G(g_{1}), \alpha_G(g_{2})  \right\}  \right)  = \left\{    q_{1},  q_{2}\right\}  \centerdot \left\{    \alpha_G(g_{1}), \alpha_G(g_{2}) \right\}  =$

 $\left[  \sigma\left\{  q_{1}, q_{2}  \right\},i \left\{   \alpha_G(g_{1}),  \alpha_G(g_{2})  \right\}  \right] = \left[  \left\{    s\left(  q_{1}\right), s\left(  q_{2}\right)   \right\}, \left\{ \alpha_G(g_{1}), \alpha_G(g_{2})  \right\}
\right]  =$

$\left\{  s\left[  q_{1},q_{2}\right]  , \alpha_G\left[g_{1},g_{2}\right]  \right\}$
\medskip

\noindent  $\rho:Ker(\pi) \otimes \frak{uce}(Q)\to Ker(\pi)$,

 $\rho\left(  \left\{   \alpha_G(g_{1}), \alpha_G(g_{2}) \right\}  \otimes\left\{  q_{1},  q_{2}  \right\}  \right)  = \left\{    \alpha_G(g_{1}), \alpha_G(g_{2}) \right\} \centerdot \left\{   q_{1},   q_{2}\right\} =$

  $\left[ i \left\{  \alpha_G(g_{1}), \alpha_G(g_{2})  \right\} ,\sigma\left\{  q_{1},  q_{2} \right\}  \right] = \left[  \left\{  \alpha_G(g_{1}), \alpha_G(g_{2})  \right\},\left\{    s\left(  q_{1}\right)  , s\left(  q_{2}\right)  \right\}  \right] =$

  $\left\{ \alpha_G \left[  g_{1},g_{2}\right]  ,s\left[q_{1},q_{2}\right]  \right\}$
\medskip

By Lemma \ref{ext rota}, the split exact sequence
 \[\xymatrix{
0 \ar[r]& (Ker(\pi), \overline{\alpha_{G }}_{\mid}) \ar[r]^{i}& (\frak{uce}_{\alpha}(G),\overline{\alpha_G})  \ar@<0.5ex>[r]^{\pi \quad} &(\frak{uce}(Q),Id_{\frak{uce}(Q)})  \ar@<0.5ex>[l]^{\sigma \quad} \ar[r]& 0} \]
is equivalent to the semi-direct product sequence, i.e.
$$\left(  \frak{uce}_{\alpha}(G),\overline{\alpha_{G}}\right)  \cong \left( Ker(\pi), \overline{\alpha_{G}}_{\mid} \right)  \rtimes \left(\frak{uce}_{\alpha}(Q),Id_{\frak{uce}_{\alpha}(Q)}\right)$$

Let $q\in Q$ and $\alpha_M(m_{1}),\alpha_M(m_{2})\in \alpha_M(M)$ be, then the following identities hold in $(\frak{uce}_{\alpha}(G), \overline{\alpha_G})$:

\noindent $\left\{  \alpha_{G}\left(  s\left(q\right)  \right)  ,\left[  t\left(  \alpha_M(m_{1})\right)  ,t\left( \alpha_M(m_{2})\right)\right]  \right\}  = \left\{  \left[  s(q),t\left( \alpha_M(m_{1})\right)  \right]  ,\alpha_{G}\left(t\left( \alpha_M(m_{2})\right)  \right)  \right\}$

\noindent $- \left\{  \left[  s(q),t( \alpha_M(m_{2}))\right]  ,\alpha_{G}\left(  t\left( \alpha_M(m_{1})\right)  \right)  \right\}$

and

\noindent$\left\{  \left[  t\left(  \alpha_M(m_{1})\right)  ,t\left(  \alpha_M(m_{2})\right)  \right]
,\alpha_{G}\left(  s\left(  q\right)  \right)  \right\}  = \left\{  \alpha_{G}\left(  t\left(  \alpha_M(m_{1})\right)  \right)  ,\left[  t\left( \alpha_M(m_{2})\right)  ,s\left(  q\right)  \right]  \right\}$

\noindent $+ \left\{  \left[ t(\alpha_M(m_{1})),s(q)\right]  ,\alpha_{G}\left(  t\left(  \alpha_M(m_{2})\right)  \right)
\right\}.$

These equalities together with the $\alpha$-perfection of $\left(  M,\alpha_{M}\right)$ imply:

\noindent $\{s(Q),M\} =\left\{  \alpha_{G}\left(  s\left(  Q\right)  \right)  , [\alpha_M(M), \alpha_M(M)] \right\} \subseteq
\left\{ \alpha_M(M), \alpha_{M}^2\left(  M\right)  \right\} \subseteq$

\noindent $\left\{\alpha_M(M),   \alpha_M(M) \right\}$

and

\noindent $\{M,s(Q)\} = \left\{  [\alpha_M(M), \alpha_M((M)], \alpha_{G}\left(  s\left(  Q\right)  \right)  \right\} \subseteq
\left\{ \alpha_M^2(M),  \alpha_M(M) \right\}+$

 \noindent $\left\{ \alpha_M(M), \alpha_{M}^2\left(  M\right)  \right\}\subseteq \left\{ \alpha_M(M), \alpha_M(M) \right\}$.

Moreover
\begin{equation} \label{Eq 1}
\tau\left(  \frak{uce}_{\alpha}(M), \overline{\alpha_M} \right) \equiv  \left( \left\{  \alpha_M(M), \alpha_M(M)  \right\}, \overline{\alpha_G}_{\mid} \right)
\end{equation}
since $\tau\left\{  \alpha_M(m_{1}) , \alpha_M(m_{2}) \right\}  = \left\{  t\left(   \alpha_M(m_{1}) \right)  ,t\left(
  \alpha_M(m_{2})\right)  \right\}  \equiv\left\{ \alpha_M(m_{1}), \alpha_M(m_{2}) \right\}$, and $$\sigma(\frak{uce}(Q))=\{s(Q),s(Q)\}=\{\alpha_G (s(Q)), \alpha_G (s(Q)) \}$$
since $\sigma(\{q_1,q_2\}) = \{s(q_1),s(q_2) \} =  \{\alpha_G(s(q_1)),\alpha_G(s(q_2)) \}$.
\medskip

On the other hand, for every $\alpha_G(g) \in G$, there exists an $\alpha_M(m) \in \alpha_M(M)$ such that $\alpha_G(g) =s \left(  p\left(    \alpha_G(g) \right)  \right)+ \alpha_M(m)$. Hence

\begin{equation} \label{Eq 2} \left( \frak{uce}_{\alpha} \left(  G\right), \overline{\alpha_G} \right)  = \left( \left\{  s\left(  Q\right),s\left(  Q\right)  \right\}  +\left\{ \alpha_M(M), \alpha_M(M) \right\}, \overline{\alpha_G} \right)
\end{equation}

\begin{Pro}\
 $$\left( Ker(\pi), \overline{\alpha_G}_{\mid} \right)= \left( \left \{ \alpha_M(M), \alpha_M(M) \right \}, \overline{\alpha_G}_{\mid} \right) = \tau\left(  \frak{uce}_{\alpha}\left(M \right),\overline{\alpha_{M}}\right).$$
\end{Pro}
{\it Proof.}  Let $\{g_1,g_2\} \in Ker(\pi)$ be. From (\ref{Eq 2}),  $\{g_1,g_2\} = \left\{  s\left(  q_{1}  \right)  ,s\left( q_{2}  \right)  \right\}  +\left\{ \alpha_M(m_{1}),\right.$  $\left. \alpha_M(m_{2})  \right\}  \in \frak{uce}_{\alpha} \left(  G\right)$. Then $\overline{0} = \pi \{g_1,g_2\} = \left\{  p\left(  s \left(  q_{1}\right) \right)  ,p\left( s\left(  q_{2}\right)  \right)  \right\} +  \left\{  p\left( \alpha_M(m_{1})  \right),p\left(  \alpha_M(m_{2})  \right)  \right\}  = \{q_1,q_2\}$, i.e. $q_{1} \otimes q_{2}   \in I_{Q}$. Consequently, $\sigma\left\{ q_{1}, q_{2}  \right\}  = \left\{   s\left(
q_{1} \right)  , s\left(  q_{2}  \right) \right\}=0$ since $s\left(  q_{1} \right) \otimes  s\left(  q_{2}  \right)  \in \sigma(I_Q) \subseteq I_{G}$. So any element in the kernel has the form  $\left\{  \alpha_M(m_{1}), \alpha_M(m_{2})  \right\}$. The reverse inclusion is obvious.

Second equality was proved in (\ref{Eq 1}). \rdg
\bigskip

On the other hand $\sigma\left(  \frak{uce}(Q), Id_{\frak{uce}(Q)} \right)  = \left( \left\{  s(Q)  ,s(Q)
  \right\}, \overline{\alpha_G} \right)$.

Since $\pi \cdot \sigma=Id_{\frak{uce}(Q)}$, then $\left( \frak{uce}_{\alpha}(G), \overline{\alpha_G} \right) = \left( Ker(\pi), \overline{\alpha_G}_{\mid} \right) \rtimes \sigma \left(  \frak{uce}\left(  Q\right), Id_{\frak{uce}(Q)} \right)$. Moreover $\sigma$ is an isomorphism between $\left( \frak{uce}(Q), Id_{\frak{uce}(Q)} \right)$ and $\sigma\left( \frak{uce}(Q), Id_{\frak{uce}(Q)} \right)$.

These facts imply:
\begin{enumerate}
\item[{\bf 1.}] $\left(  \frak{uce}_{\alpha}(G),\overline{\alpha_{G}} \right)  =\tau\left(
\frak{uce}_{\alpha} \left(  M\right), \overline{\alpha_{M}} \right)  \rtimes
\sigma\left(  \frak{uce}\left(  Q\right), Id_{\frak{uce}(Q)}\right).$

\item[{\bf 2.}]  $\sigma\left(  \frak{uce} \left(  Q\right),Id_{\frak{uce}(Q)} \right)  \cong\left(  \frak{uce} \left(  Q\right), Id_{\frak{uce}(Q)}\right).$
\end{enumerate}

From {\bf 1.}, an element of $\left(  \frak{uce}_{\alpha}(G),\overline{\alpha_{G}} \right)$ can be written as
$\left(  \tau\left(  m\right)  ,\sigma\left(  q\right)  \right)$, for
$m\in\left(  \frak{uce}_{\alpha}\left(  M\right),\overline{\alpha_{M}}\right)$ and
$q\in\left(  \frak{uce}\left(  Q\right), Id_{\frak{uce}(Q)} \right)$ with a suitable choice.
 Such an element belongs to $Ker(U_{\alpha}^G)$ if and only if $U_{G}^{\alpha}\left(  \tau\left(  m\right)  ,\sigma\left(  q\right)  \right)  =0$, i.e. $m\in Ker(U_{\alpha}^M)$ and $q\in HL_2(Q)$.

From these facts we can derive that
\begin{enumerate}
\item[{\bf 3.}] $\left( Ker(U_{\alpha}^G), \overline{\alpha_G}_{\mid} \right) \cong \tau\left( Ker(U_{\alpha}^M) , \overline{\alpha_{M}}_{\mid} \right)  \oplus \sigma \left(  HL_2(Q), {Id_{\frak{uce}(Q)}}_{\mid}  \right).$
\end{enumerate}

Since there exists a symmetric Hom-action of $(Q,Id_Q)$ on $(M, \alpha_M)$, then there is a Hom-action of
$\left(  \frak{uce}(Q), Id_{\frak{uce}(Q)} \right)$ on $\left(  \frak{uce}_{\alpha}(M),\overline{\alpha_{M}}\right)$ given by:

$\begin{array}{rcl} \lambda : \frak{uce}(Q) \otimes \frak{uce}_{\alpha}(M)& \to &  \frak{uce}_{\alpha}(M)\\
  \left\{  q_{1}, q_{2} \right\}  \otimes \left\{ \alpha_M(m_{1}), \alpha_M(m_{2}) \right\}   & \mapsto & \left\{ q_{1},  q_{2}\right\}  \centerdot \left\{  \alpha_M(m_1), \alpha_M(m_{2})  \right\}  =\\
 & & \left\{  \left[  q_{1},q_{2}\right]  \centerdot  \alpha_M(m_{1}),\alpha_{M}^2 \left(  m_{2}\right)  \right\}  -\\
 & & \left\{  \left[  q_{1},q_{2}\right]  \centerdot  \alpha_M(m_{2}),\alpha_{M}^2\left(  m_{1}\right)  \right\}
\end{array}$

and

$\begin{array}{rcl} \rho: \frak{uce}_{\alpha}(M) \otimes \frak{uce}(Q) & \to & \frak{uce}_{\alpha}(M)\\
 \left\{ \alpha_M(m_{1}), \alpha_M(m_{2})  \right\}  \otimes \left\{   q_{1}, q_{2} \right\}  &\mapsto & \left\{    \alpha_M(m_{1}),  \alpha_M(m_{2})\right\} \centerdot \left\{ q_{1}, q_{2}  \right\} = \\
& & \left\{  \alpha_M(m_{1})  \centerdot \left[  q_{1},q_{2}\right],\alpha_{M}^2\left(  m_{2}\right)  \right\} -\\
&&\left\{  \alpha_{M}^2\left(  m_{1}\right)  ,\left[  q_{1},q_{2}\right] \centerdot  \alpha_M(m_{2})  \right\}
\end{array}$
\medskip

Then we can define the following homomorphism of  Hom-Leibniz algebras:
$$\tau \rtimes \sigma : \left(  \frak{uce}_{\alpha}\left(  M\right), \overline{\alpha_{M}}\right)  \rtimes \left(  \frak{uce} \left(  Q\right), Id_{\frak{uce}(Q)}\right) \to \left(  \frak{uce}_{\alpha}\left(  G\right) , \overline {\alpha_{G}}\right)  \cong \quad \quad \quad \quad \quad $$
$$\quad \quad \quad \quad \quad \quad \quad \quad \quad \quad \quad \quad \quad \quad \quad \quad \quad \quad \tau \left(  \frak{uce}_{\alpha}\left(  M\right),  \overline{\alpha_{M}}\right)  \rtimes \sigma \left(  \frak{uce}\left( Q\right), Id_{\frak{uce}(Q)}\right)$$
$$\left(  \left\{  \alpha_M(m_{1}), \alpha_M(m_{2})  \right\},\left\{ q_{1},  q_{2} \right\}  \right) \mapsto \left(  \left\{
t\left( \alpha_M(m_{1}) \right)  ,t\left(  \alpha_M(m_{2}) \right)  \right\},\left\{  s\left(  q_{1}  \right),s\left(    q_{2} \right) \right\}  \right)$$

 Moreover $\tau\rtimes\sigma$ is an epimorphisms since $$\left( \frak{uce}_{\alpha}(G),\overline{\alpha_{G}}\right)  \cong\tau\left(  \frak{uce}_{\alpha}\left(  M\right),\overline{\alpha_{M}}\right)  \rtimes\sigma\left(  \frak{uce}\left(  Q\right),Id_{\frak{uce}(Q)}\right).$$

By the relations coming from the action induced by the  split extension

\noindent $\tau\left(  \left\{  q_{1},q_{2}\right\}  \centerdot \left\{\alpha_M(m_{1}), \alpha_M(m_{2})\right\}  \right)  =\left[  \left\{  s\left(  q_{1}\right)
,s\left(  q_{2}\right)  \right\}, \left\{ t\left(  \alpha_M(m_{1})\right)  ,t\left( \alpha_M(m_{2})\right) \right\}  \right]$

and

\noindent $\tau\left(  \left\{  \alpha_M(m_{1}), \alpha_M(m_{2})\right\}  \centerdot \left\{  q_{1},q_{2}\right\}
\right)  =$ $\left[  \left\{  t\left(  \alpha_M(m_{1})\right), t\left(  \alpha_M(m_{2})\right)
\right\}  ,\left\{  s\left(  q_{1}\right)  ,s\left(  q_{2}\right)  \right\}
\right]$

one derives that:

\noindent $t \cdot U_{\alpha}^{M}\left(  \left\{  q_{1},q_{2}\right\}  \centerdot \left\{  \alpha_M(m_{1}), \alpha_M(m_{2}) \right\}  \right)=
 \left[q_1,q_2 \right] \centerdot [\alpha_M(m_1), \alpha_M(m_2)],$

and

\noindent $t \cdot U_{\alpha}^{M}\left(  \left\{  \alpha_M(m_{1}), \alpha_M(m_{2}) \right\}  \centerdot \left\{  q_{1},q_{2}\right\}  \right)  =
 \left[\alpha_M(m_1), \alpha_M(m_2) \right] \centerdot [q_1,q_2].$

\medskip

\begin{enumerate}
\item[{\bf 4.}]
Now we define the surjective homomorphism of Hom-Leibniz algebras
\end{enumerate}
\begin{center}
$\begin{array}{rcl}
\Phi:= \left(  t \cdot U_{\alpha}^{M}\right)  \rtimes\left(  s \cdot u_{Q}\right):\left(
\frak{uce}_{\alpha} \left(  M\right)  \rtimes \frak{uce} \left(  Q\right)  ,\overline{\alpha_{M}}\rtimes Id_{\frak{uce}(Q)}\right) & \to & \left(  G,\alpha_{G}\right)  \\
 \left(  \left\{  \alpha_M(m_{1}), \alpha_M(m_{2}) \right\},\left\{  q_{1},q_{2}\right\}  \right)  & \mapsto & \left(  t\left[
\alpha_M(m_{1}), \alpha_M(m_{2})\right]  ,s\left[  q_{1},q_{2}\right]  \right) \end{array}$
\end{center}

that makes commutative the following diagram:
\begin{equation} \label{diagrama Fi}
\vcenter{ \xymatrix{
\left( \frak{uce}_{\alpha} \left(  M\right)  \rtimes \frak{uce} \left(  Q\right),\overline{\alpha_{M}}\rtimes Id_{\frak{uce}(Q)}\right)  \ar@{>>}[rr]^{\quad \tau \rtimes \sigma} \ar@{-->}[dr]_{\Phi}& &  \left(  \frak{uce}_{\alpha} \left(  G\right),\overline{\alpha_{G}} \right)\ar@{>>}[ld]^{U_{\alpha}^{G}}\\
&  \left(  G,\alpha_{G}\right)&
}}
\end{equation}

Now we prove that $$\frak{uce} \left(  Q\right) \centerdot Ker(U_{\alpha}^M)  \oplus  Ker(U_{\alpha}^M)  \centerdot \frak{uce}\left(  Q\right)   \subseteq Ker(\tau) \subseteq Ker(U_{\alpha}^M)$$

Second inclusion is obvious since $t \cdot U_{\alpha}^M = U_{\alpha}^G \cdot \tau$ and $t$ is injective.

From the commutativity of the following diagram
\[ \xymatrix{
 & \left(  Ker(U_{\alpha}^M),\overline{\alpha_{M}}_{\mid}\right)   \ar@{-->}[r]  \ar@{>->}[d] & \left(  Ker(U_{\alpha}^G), \overline{\alpha_{G}}_{\mid}\right) \ar@{>->}[d]\\
Ker(\tau)\ \ar@{-->}[ru] \ar@{>->}[r] & \left( \frak{uce}_{\alpha}(M),\overline{\alpha_{M}} \right) \ar[r]^{\tau} \ar@{>>}[d]^{U_{\alpha}^M}& \left(  \frak{uce}_{\alpha}(G),\overline{\alpha_{G}} \right) \ar@{>>}[d]^{U_{\alpha}^G}\\
  & \left(  M,\alpha_{M}\right)  \ar[r]^t & \left(G,\alpha_{G}\right)
} \]
we have that $U_{\alpha}^{G} \cdot \tau\left(  Ker(U_{\alpha}^M) \right)  = t \cdot U_{\alpha}^{M}\left( Ker(U_{\alpha}^M)\right)
=0$, then $\tau\left(  Ker(U_{\alpha}^M) \right)$ $\subseteq Ker(U_{\alpha}^{G}) \subseteq Z\left(  \frak{uce}_{\alpha}(G)\right)$, so,

$\tau\left( \frak{uce}\left(  Q\right)  \centerdot Ker(U_{\alpha}^M) \right)  =\left[  \sigma\left(  \frak{uce} \left(  Q\right)  \right)  ,\tau\left( Ker(U_{\alpha}^M)  \right)  \right]  =0$

\noindent and

$\tau\left(  Ker(U_{\alpha}^M)  \centerdot  \frak{uce}\left(  Q\right)\right)  =\left[  \tau\left(  Ker(U_{\alpha}^M) \right),\sigma\left(  \frak{uce}\left(  Q\right)  \right)  \right]  =0$

Consequently, $\frak{uce} \left(  Q\right)  \centerdot  Ker(U_{\alpha}^M)  \oplus  Ker(U_{\alpha}^M)  \centerdot  \frak{uce}\left(  Q\right) \subseteq Ker(\tau).$
\medskip

On the other hand, we observe that  $\left( \frak{uce}\left(  Q\right)  \centerdot Ker(U_{\alpha}^M)  \oplus   Ker(U_{\alpha}^M)  \centerdot \frak{uce}\left(  Q\right), \right.$ $\left. \overline{\alpha_M}_{\mid} \right)$ is  a two-sided  ideal of $\left( \frak{uce}_{\alpha}\left(  M\right), \overline{\alpha_M} \right)$. Then the Hom-action of $\left( \frak{uce}\left(  Q\right), Id_Q \right)$ on $\left( \frak{uce}_{\alpha}\left(  M\right), \overline{\alpha_M} \right)$ induces a Hom-action of $\left( \frak{uce}\left(  Q\right), Id_Q \right)$ on $$\left( \overline{\frak{uce}_{\alpha}\left(  M\right)  }, \overline{\overline{\alpha_M}} \right)= \left( \frac{\frak{uce}_{\alpha}\left(  M\right)}{\frak{uce}\left(  Q\right)  \centerdot  Ker(U_{\alpha}^{M})
\oplus  Ker(U_{\alpha}^{M}) \centerdot \frak{uce}\left(  Q\right)}, \overline{\overline{\alpha_M}} \right).$$

Since $\tau$ vanishes on $\frak{uce}\left(  Q\right)  \centerdot Ker(U_{\alpha}^{M}) \oplus  Ker(U_{\alpha}^{M}) \centerdot
\frak{uce} \left(  Q\right)$, then it induces $\overline{\tau}:\overline{\frak{uce}_{\alpha}\left(  M\right)}\to
\tau\left(  \frak{uce}_{\alpha}\left(  M\right)  \right)$. This fact is illustrated in the following diagram where the notation
 $I=  \frak{uce} \left(  Q\right)  \centerdot Ker(U_{\alpha}^{M})  \oplus Ker(U_{\alpha}^{M}) \centerdot  \frak{uce}\left(Q\right)$ is employed:
\[ \xymatrix{
\left( I,  \overline{\alpha_M}_{\mid} \right) \ar@{>->}[d] \ar@/^1pc/[ddr]^0&  &\\
\left(  \frak{uce}_{\alpha}(M),\overline{\alpha_{M}}\right)  \ar[rr]^{\tau} \ar@{>>}[dd] \ar@{>>}[dr]& & \left( \frak{uce}_{\alpha}(G),\overline{\alpha_{G}}\right) \\
& \tau \left( \frak{uce}_{\alpha}(M), \overline{\alpha_M} \right) \ar@{>->}[ur]&\\
\left( \overline{\frak{uce}_{\alpha}\left(  M\right)}, \overline{\overline{\alpha_M}} \right)\ar@{-->}[ur]^{\overline{\tau}} & &
} \]
Now we can construct the following commutative diagram:
\[ \xymatrix{
I\  \ar@{>->}[r] \ar@{>->}[d] & I\rtimes 0  \ar@{>>}[r] \ar@{>->}[d] & 0 \ar@{>->}[d] \\
Ker\left(  \tau\rtimes\sigma\right)\ \ar@{>>}[d] \ar@{>->}[r] & \left(  \frak{uce}_{\alpha
}\left(  M\right)  \rtimes \frak{uce}\left(  Q\right)  ,\overline{\alpha
_{M}}\rtimes Id_{\frak{uce}(Q)}\right)  \ar@{>>}[r]^{\quad \quad \quad \quad \quad \tau\rtimes\sigma} \ar@{>>}[d]& \left(  \frak{uce}_{\alpha}\left(  G\right)  ,\overline{\alpha_{G}}\right) \ar@{=}[d]\\
\frac{Ker\left(  \tau\rtimes\sigma\right)}{I}\ \ar@{>->}[r] & \left(  \overline{\frak{uce}_{\alpha}\left(  M\right)  }\rtimes \frak{uce}\left(  Q\right),\overline{\overline{\alpha_{M}}}\rtimes Id_{\frak{uce}(Q)}\right)  \ar@{>>}[r]^{\quad \quad \quad \quad \quad \Psi} & \left(  \frak{uce}_{\alpha}\left(  G\right)  ,\overline{\alpha_{G}}\right)
}\]
whose bottom row is a central extension. Moreover  $\left(  \frak{uce}_{\alpha}\left(  G\right),\overline{\alpha_{G}}\right)$ is
an $\alpha$-perfect Hom-Leibniz algebra, then by  Theorem 5.5 in \cite{CIP}, it admits a universal $\alpha$-central extension and, by Corollary \ref{centralmente cerrada},  $\frak{uce}_{\alpha}\left(  G\right)$ is centrally closed,  i.e. $\frak{uce} \left(  \frak{uce}_{\alpha}\left(  G\right)  \right)  \cong \frak{uce}_{\alpha}\left(  G\right)$.

Having in mind the following diagram,
\[ \xymatrix{
\left(  \overline{\frak{uce}_{\alpha}\left(  M\right)  }\rtimes \frak{uce}\left( Q\right),\overline{\alpha_{M}}\rtimes Id_{\frak{uce}(Q)}\right) \ar@{>>}[r]^{\quad \quad \quad \quad   \Psi} \ar@{>>}[d]^{\Psi}  \ar@/_4pc/[dd]_{Id} & \left(  \frak{uce}_{\alpha}\left(  G\right),\overline{\alpha_{G}}\right) \ar@{=}[d]\\
\left(  \frak{uce}_{\alpha} \left(  G\right)  ,\overline{\alpha_{G}}\right)  \ar@{>->>}[r]^{\quad  Id} \ar[d]^{\mu} & \left(  \frak{uce}_{\alpha}\left(  G\right)  ,\overline{\alpha_{G}}\right) \ar@{=}[d] \\
\left(  \overline{\frak{uce}_{\alpha}\left(  M\right)  }\rtimes \frak{uce}\left(Q\right)  ,\overline{\alpha_{M}}\rtimes Id_{\frak{uce}(Q)}\right) \ar@{>>}[r]^{\quad \quad \quad \quad \Psi} & \left(  \frak{uce}_{\alpha}\left(  G\right)
,\overline{\alpha_{G}}\right)
} \]
where $Id:\left(  \frak{uce}_{\alpha}\left(  G\right),\overline{\alpha_{G}}\right)  \to \left(  \frak{uce}_{\alpha}\left(
G\right)  ,\overline{\alpha_{G}}\right)$ is a universal central extension since $\left(  \frak{uce}_{\alpha}\left(  G\right)  ,\overline{\alpha_{G}}\right)$ is centrally closed and $\Psi$ is a central extension, then there exists a unique
 homomorphism of Hom-Leibniz algebras $\mu: \left( \frak{uce}_{\alpha}\left( G\right), \overline{\alpha_G} \right)  \to \left( \overline{\frak{uce}_{\alpha}\left(  M\right)  }\rtimes \frak{uce}\left(  Q\right), \overline{\alpha_{M}}\rtimes Id_{\frak{uce}(Q)}\right)$ such that $\Psi \cdot \mu=Id.$

Since $\Psi \cdot \mu \cdot \Psi=Id \cdot \Psi=\Psi=\Psi \cdot Id$ and $\overline{\frak{uce}_{\alpha} \left(  M\right)  }\rtimes \frak{uce}\left(  Q\right)$ is $\alpha$-perfect, then Lemma 5.4 in \cite{CIP} implies that $\mu \cdot \Psi=Id$.
Consequently, $\Psi$ is an isomorphism, then $Ker(\Psi)= \frac{Ker\left(  \tau\rtimes\sigma\right)}{I}=0$, so $Ker\left(  \tau\rtimes\sigma\right)  \subseteq I$.

The above discussion can be summarized in:

\begin{enumerate}
\item[{\bf 5.}] $Ker\left(  \tau\rtimes\sigma\right)  \cong   \frak{uce}\left(  Q\right)  \centerdot  Ker(U_{\alpha}^{M})  \oplus Ker(U_{\alpha}^{M}) \centerdot \frak{uce}\left(  Q\right)$
\end{enumerate}

\medskip

We summarize the above results in the following

\begin{Th} \label{5 puntos}
Consider a split extension of $\alpha$-perfect Hom-Leibniz algebras
\[\xymatrix{
0 \ar[r]& (M,\alpha_M) \ar[r]^t& (G,\alpha_G) \ar@<0.5ex>[r]^p &(Q,Id_Q) \ar[r] \ar@<0.5ex>[l]^s &0} \]
where the induced Hom-action of $(Q, Id_Q)$ on $(M, \alpha_M)$ is symmetric. Then the following statements hold:
\begin{enumerate}
\item[{\bf 1.}] $\left(  \frak{uce}_{\alpha}(G),\overline{\alpha_{G}} \right)  =\tau\left(
\frak{uce}_{\alpha} \left(  M\right), \overline{\alpha_{M}} \right)  \rtimes
\sigma\left(  \frak{uce}\left(  Q\right), Id_{\frak{uce}(Q)}\right).$

\item[{\bf 2.}]  $\sigma\left(  \frak{uce} \left(  Q\right),Id_{\frak{uce}(Q)} \right)  \cong\left(  \frak{uce} \left(  Q\right), Id_{\frak{uce}(Q)}\right).$

    \item[{\bf 3.}] $\left( Ker(U_{\alpha}^G), \overline{\alpha_G}_{\mid} \right) \cong \tau\left( Ker(U_{\alpha}^M) , \overline{\alpha_{M}}_{\mid} \right)  \oplus \sigma \left(  HL_2(Q), {Id_{\frak{uce}(Q)}}_{\mid}  \right).$

 \item[{\bf 4.}]
The homomorphism of Hom-Leibniz algebras $$\Phi:\left(\frak{uce}_{\alpha} \left(  M\right)  \rtimes \frak{uce} \left(  Q\right), \overline{\alpha_{M}}\rtimes Id_{\frak{uce}(Q)}\right)\to  \left(  G,\alpha_{G}\right)$$ given by $\Phi
 \left(  \left\{  \alpha_M(m_{1}), \alpha_M(m_{2}) \right\},\left\{  q_{1},q_{2}\right\}  \right)  = \left(  t\left[
\alpha_M(m_{1}), \alpha_M(m_{2})\right]  ,s\left[  q_{1},q_{2}\right]  \right)$ is an epimorphism that makes commutative diagram (\ref{diagrama Fi}) and its kernel is $Ker(U_{\alpha}^{M})\oplus HL_2(Q)$.

   \item[{\bf 5.}] $Ker\left(  \tau\rtimes\sigma\right)  \cong   \frak{uce}\left(  Q\right)  \centerdot Ker(U_{\alpha}^{M})  \oplus Ker(U_{\alpha}^{M}) \centerdot \frak{uce}\left(  Q\right)$
\end{enumerate}
\end{Th}

\begin{Rem}
Let us observe that statements {\bf 1., 2.} and {\bf 3.} in Theorem \ref{5 puntos} hold in the general case, they do not need the hypothesis of symmetric Hom-action.
\end{Rem}

\begin{Th} The following statements are equivalent:
\begin{enumerate}
\item[a)]  $\Phi = (t \cdot U_{\alpha}^M) \rtimes (s \cdot u_Q) :\left(  \frak{uce}_{\alpha}\left(  M\right)  \rtimes
\frak{uce}\left(  Q\right)  ,\overline{\alpha_{M}}\rtimes Id_{\frak{uce}_{\alpha}(Q)}\right)  \to \left(  G,\alpha_{G}\right)$ is a central extension, hence is an $\alpha$-cover.

\item[b)] The Hom-action of $(\frak{uce}(Q), Id_Q)$ on $(Ker(U_{_{\alpha}}^{M}), \overline{\alpha_M}_{\mid})$ is trivial.

\item[c)]  $\tau \rtimes \sigma$ is an isomorphism. Consequently $\frak{uce}_{\alpha}\left(  M\right)  \rtimes \frak{uce}\left(  Q\right)$  is the universal $\alpha$-central extension of $(G, \alpha_G)$.

\item[d)] $\tau$ is injective.
\end{enumerate}

 In particular, for the direct product $\left(  G,\alpha_{G}\right)  = \left(  M,\alpha_{M} \right) \times \left(Q, Id_{Q}\right)$ the following isomorphism hold:
$$\left(  \frak{uce}_{\alpha}\left(  M\times Q\right),\overline{\alpha_{M}\times Id_{Q}}\right)  \cong \left(  \frak{uce}_{\alpha}\left(  M\right)  \times \frak{uce}\left(  Q\right) ,\overline{\alpha_{M}}\times Id_{\frak{uce}(Q)}\right)$$
\end{Th}
{\it Proof.}
{\it a)} $\Longleftrightarrow$ {\it b)}

 If $\Phi:\left(  \frak{uce}_{\alpha}\left(  M\right)  \rtimes \frak{uce}_{\alpha}\left(  Q\right)  ,\overline{\alpha_{M}}\rtimes Id_{\frak{uce}_{\alpha}(Q)}\right) \to \left(  G,\alpha_{G}\right)$ is a central extension and having in mind that  $Ker(\Phi) = Ker(U_{\alpha}^{M})\oplus HL_2(Q)$, then the Hom-action of $(\frak{uce}(Q), Id_Q)$ on $(Ker(U_{\alpha}^M), \overline{\alpha_M}_{\mid})$ is trivial and conversely.

Moreover $(\frak{uce}_{\alpha}\left(  M\right)  \rtimes \frak{uce}\left(Q\right), \overline{\alpha_M} \rtimes Id_{\frak{uce}(Q)})$ is $\alpha-$perfect since the Hom-action is trivial.

{\it b)} $\Longleftrightarrow$ {\it c)}

By statement {\bf 5.} in Theorem \ref{5 puntos}   we know that  $Ker\left(  \tau\rtimes\sigma\right)  \cong \frak{uce}\left(  Q\right) \centerdot Ker(U_{_{\alpha}}^{M}) \oplus  Ker(U_{_{\alpha}}^{M}) \centerdot \frak{uce}\left(  Q\right)$, then $\tau \rtimes \sigma$ is injective if and only if the Hom-action is trivial.

Hence and having in mind  diagram (\ref{diagrama Fi}), immediately follows that  $\frak{uce}_{\alpha}\left(  M\right)  \rtimes \frak{uce}\left(  Q\right)$ is the universal $\alpha$-central extension of $(G, \alpha_G)$.

{\it c)} $\Longleftrightarrow$ {\it d)}

It suffices to have in mind  the identification  of $\tau$ with $\tau \rtimes \sigma$ given by $\tau\left\{  \alpha_{M}\left(  m_{1}\right)  ,\right.$ $\left.
\alpha_{M}\left(  m_{2}\right)  \right\}  =\left( \tau\rtimes\sigma\right)  \left(  \left\{  \alpha_{M}\left(  m_{1}\right),\alpha_{M}\left(  m_{2}\right)  \right\}  ,0\right)$, since $Ker(\tau) \cong Ker\left(  \tau\rtimes\sigma\right)$, then the equivalence is obvious.
\medskip

Finally, in case of the direct product $\left(  G,\alpha_{G}\right)  =\left(M, \alpha_{M}\right) \times \left(Q, Id_{Q}\right)$ the Hom-action of $\left(Q,Id_{Q}\right)$ on $\left(  M,\alpha_{M}\right)$ is trivial, then the Hom-action of $(\frak{uce}(Q), Id_{\frak{uce}(Q)})$ on $(\frak{uce}_{\alpha}(M), \overline{\alpha_M})$  is trivial as well and, consequently,  $(\frak{uce}_{\alpha}(M)\rtimes \frak{uce}(Q), \overline{\alpha_M} \rtimes Id_{\frak{uce}(Q)})=(\frak{uce}_{\alpha}(M)\rtimes \frak{uce}(Q), \overline{\alpha_M} \times Id_{\frak{uce}(Q)})$.

Statement  {\it c)} ends the proof. \rdg

\begin{Rem}
Note that when the Hom-Leibniz algebras are considered as  Leibniz algebras, i.e. the endomorphisms $\alpha$ are identities, then the results in this section recover the  corresponding results for Leibniz algebras given in \cite{CC}.
\end{Rem}

\begin{center}

\end{center}

\end{document}